\documentclass[11pt,a4paper]{article}
\usepackage{format}

\title{Limit Theorems for One-Dimensional Homogenized Diffusion Processes\thanks{Submitted to \emph{Journal of Statistical Physics} on October 22, 2025.}}
\author{Jaroslav I. Borodavka, Sebastian Krumscheid}
\date{October 22, 2025}

\begin{document}

\maketitle

\begin{abstract}
\noindent
We present two limit theorems, a mean ergodic and a central limit theorem, for a specific class of one-dimensional diffusion processes that depend on a small-scale parameter $\eps$ and converge weakly to a homogenized diffusion process in the limit $\epstozero$. In these results, we allow for the time horizon to blow up such that $T_\eps \rightarrow \infty$ as $\epstozero$. The novelty of the results arises from the circumstance that many quantities are unbounded for $\epstozero$, so that formerly established theory is not directly applicable here and a careful investigation of all relevant $\eps$-dependent terms is required. As a mathematical application, we then use these limit theorems to prove asymptotic properties of a minimum distance estimator for parameters in a homogenized diffusion equation.
\end{abstract}

\section{Introduction} \label{sec:introduction}

Dating back to Boltzmann's ergodic hypothesis essentially stating that time averages converge to ensemble averages \cite{K:1949}, one of the driving principles of statistical physics is arguably the one of ergodicity. This property can be observed and postulated for many dynamical systems. Among these are also Langevin diffusions which play a central role in the field of molecular dynamics, see for example \cite{LS:2016, S:2010}, where macroscopic properties and observables are inferred from atomistic models usually under the assumption of ergodicity. In molecular dynamics and many other fields, dynamics can be characterized by processes occurring across multiple time scales. A mathematically convenient description for such dynamics is often given by stochastic differential equations (SDEs). In particular, in the seminal paper \cite{PS:2007}, the authors considered a class of multiscale overdamped Langevin diffusions as their data-generating model given by the SDE
\begin{equation}     \label{eq:multiscale_langevin}
    dX_\eps(t) = -\alpha V'(X_\eps(t)) - \frac{1}{\eps} p'\left( \frac{X_\eps(t)}{\eps} \right) dt  + \sqrt{2 \sigma} dW(t), \quad t \in [0,T],
\end{equation}
where $V$ is a large-scale potential, $p$ is a 1-periodic function, and $\eps>0$ is a parameter measuring the scale separation. This sequence of diffusion processes has a well-defined homogenization limit for $\epstozero$ given by a Langevin diffusion with new coefficients damped by a homogenization factor $K \leq 1$
\begin{equation}    \label{eq:limit_langevin}
    dX(t) = -\alpha K V'(X(t)) dt  + \sqrt{2 \sigma K} dW(t),  \quad t \in [0,T].
\end{equation}
In that paper, the authors analyzed the maximum likelihood estimation built on basis of the homogenized model equation \eqref{eq:limit_langevin} while subject to multiscale observations from \eqref{eq:multiscale_langevin} and under the aspect of subsampling. Their established limit results are sequential in nature, that is, they first let the time horizon $\Ttoinfty$ and then $\epstozero$, where one of the main theoretical tools applied here are limit theorems for diffusion processes. More precisely, if we consider the solution $X$ of the general one-dimensional Itô SDE
\begin{equation}     \label{eq:original_sde_intro}
    d X(t) = b(X(t)) dt + \overline{\sigma}(X(t)) d\overline{W}(t), \quad t \in [0, T],
\end{equation}
and suppose $X$ has an invariant density $\mu$ on $\R$, then, under various assumptions, there are limit theorems in the existing literature, cf. \cite{K:2004, DZ:1996}, of the type
\begin{equation}
    \frac{1}{T} \int_0^T h(X(t)) \, dt \rightarrow \int_\R h(x) \mu(x) \, dx, \quad T \rightarrow \infty,
\end{equation}
for a suitable test function $h$ with the convergence understood either in the mean or almost sure sense, or as convergence in probability. Another closely related limit result is the central limit theorem (CLT)
\begin{equation}
    \frac{1}{\sqrt{T}} \int_0^T h(X(t)) \, dt \weakconv \mathcal{N}(0, \tau^2), \quad T \rightarrow \infty,
\end{equation}
for a function $h$ centered with respect to $\mu$, where $\mathcal{N}(0, \tau^2)$ is a centered normal distribution with variance $\tau^2>0$. Our main contribution in this paper is to extend the preceding two types of limit theorems in the following way for a particular class of diffusion processes, including \eqref{eq:multiscale_langevin}, given by the solution $X_\eps$ of the one-dimensional SDE
\begin{equation}    \label{eq:original_sde_eps_intro}
    d X_\eps(t) = \left[ f_0(X_\eps(t)) + \frac{1}{\eps} f_1\left( \frac{X_\eps(t)}{\eps} \right) \right] dt + \sigma(X_\eps(t)) \, dW(t),  \quad t \in [0, T],
\end{equation}
where $\eps > 0$.\label{rev:new_references} Following the groundbreaking work of \cite{K:1966a, K:1966b, PSV:1976, BLP:1978}, the later works of \cite{P:1999, PV:2001, PS:2008}, and possibly the most recent contribution in \cite{RX:2021}, this puts us, under appropriate assumptions on the coefficient functions, in the homogenization setting for SDEs. Homogenization for SDEs is a type of diffusion approximation in which the process $X_\eps$ converges weakly in $C([0, T]; \R)$ to the process $X$ as $\epstozero$ for fixed $T>0$, where $X$ is the solution to \eqref{eq:original_sde_intro}. The homogenized coefficient functions $b$ and $\overline{\sigma}$ in \eqref{eq:original_sde_intro} depend on $f_0, f_1, \sigma$, and the invariant measure of the fast process defined through $Y_\eps := X_\eps/\eps$. In this framework, our aim is to prove the following two limit theorems. Firstly, a mean ergodic theorem (MET) of the type
\begin{equation}    \label{eq:intro_met}
    \frac{1}{T_\eps} \int_0^{T_\eps} h_\eps(X_\eps(t)) \, dt \rightarrow \int_\R h(x) \mu(x) \, dx, \quad \epstozero,
\end{equation}
for a suitable sequence of functions $h_\eps$ depending on $\eps$ such that $h_\eps \rightarrow h$, and exploding time horizon such that $T_\eps$ explicitly depends on $\eps$ with $T_\eps \rightarrow \infty$ as $\epstozero$. Secondly, we want to prove a CLT of the type
\begin{equation}     \label{eq:intro_clt}
    \frac{1}{\sqrt{T_\eps}} \int_0^{T_\eps} h_\eps(X_\eps(t)) \, dt \weakconv \mathcal{N}(0, \tau^2), \quad \epstozero.
\end{equation}
The interest in these two limit theorems stems from the increasing emergence of parametric and nonparametric estimation procedures for homogenized equations like \eqref{eq:original_sde_intro}, but with observations coming from perturbed equations like \eqref{eq:original_sde_eps_intro}. Although the laws of the processes are similar in a weak sense when $\epstozero$, it turns out that many estimators cannot recover homogenized quantities from perturbed observations, presumably due to a discontinuity with respect to the weak topology. This has been repeatedly demonstrated; cf. \cite{PS:2007, PPS:2009, MP:2018, PRZ:2025}. In the aforementioned references, as well as in \cite{K:2018, AGPSZ:2021}, the authors utilize exclusively sequential limits to prove asymptotic properties of their estimators, that is, they first let $T \rightarrow \infty$ and then $\epstozero$, or vice versa. To the best of our knowledge, the case of taking coupled limits as in \eqref{eq:intro_met} and \eqref{eq:intro_clt} where almost every relevant quantity depends on $\eps$ has not been considered yet.\label{rev:motivation} At first glance this coupling of time with the scale parameter might seem unmotivated, but, in fact, it is particularly important when rigorously proving asymptotic normality of estimators for misspecified models. To be more precise, in the statistical inference with misspecification between data-generating model and model of interest, one often has to deal with quantities like
\begin{equation*}
    \frac{1}{\sqrt{T}} \int_0^T h(X_\eps(t)) - \langle h \rangle_\mu \, dt, \quad \langle h \rangle_\mu := \int_\R h(x) \mu(x) \, dx,
\end{equation*}
but, in order to analyze this quantity for a CLT, the function $h$ needs to be centered with respect to the correct invariant density, which is $\mu_\eps$, the one associated with the process $X_\eps$. Doing so, we obtain the splitting
\begin{equation*}
    \frac{1}{\sqrt{T}} \int_0^T h(X_\eps(t)) - \langle h \rangle_\mu\, dt = \frac{1}{\sqrt{T}} \int_0^T h(X_\eps(t)) - \langle h \rangle_{\mu_\eps} \, dt + \sqrt{T} \left( \langle h \rangle_{\mu_\eps} - \langle h \rangle_\mu \right)\;.
\end{equation*}
But here comes the crux. We cannot rigorously justify letting $\Ttoinfty$ and then $\epstozero$ as in the previously mentioned references due to the second summand. We let $T$ depend on $\eps$, so that we can leverage the two factors in the second summand and obtain a rigorous and logically sound CLT. Some authors, such as in \cite{AGPSZ:2021}, derived a CLT for an estimator formally, but not on the general mathematical level as we aim for. The avoidance of analyzing coupled limits can possibly be attributed to the difficulty that many relevant quantities and bounds explode for $\epstozero$, so that classically established results in the literature for diffusion processes cannot be applied directly. In addition to that, many results concerned with the homogenization of SDEs are obtained for fixed, finite time horizons $T$ and $\epstozero$, so that existing bounds and estimates in that field of study will be rather crude and insufficient in terms of $T$ for our purposes. Hence, with this work, we hope to lay the theoretical groundwork for the rigorous asymptotic analysis of estimators for homogenized models when confronted with perturbed data.

The main contributions of our work are the proofs of the two limit theorems, a convergence in probability result stringently connected to the proof of the CLT, and a simple application to parameter estimation in a continuous-time setting for homogenized models under perturbed data. At this point, it is important to emphasize that we use different proof ideas and methods from different references, and we will take great care in crediting the authors of the original proofs; however, the $\eps$-dependence in all relevant quantities demands nonetheless a careful investigation and logical extension of existing proofs where no such $\eps$-dependence is present. That one has to practice caution with tracking relevant quantities in $\eps$ is mostly due to the fact that many of these quantities, e.g., the drift of $X_\eps$, are not bounded in $\eps$ so that, for example, even standard dissipativity assumptions fail at the $\eps$-level for $\epstozero$.

The main content of the article starts with section \ref{sec:limit_theorems}, which is itself divided into three subsections. The first subsection deals with the precise definition of the framework, the assumptions under which we will work, and the rigorous proof of the MET, which is stated in Theorem \ref{theo:eps_mean_ergodic_theorem}. Subsection \ref{subsec:prob_convergence} is devoted to the proof of a required convergence in probability result contained in Corollary \ref{cor:prob_conv_result}. The last subsection \ref{subsec:central_limit_theorem} presents the proof of the CLT whose statement can be found in Theorem \ref{theo:eps_clt}. Finally, an exemplary application to parameter estimation is discussed in section \ref{sec:parameter_estimation}, with the main result being Proposition \ref{prop:mde_exa_result}. Finally, in the conclusion, we summarize the main contributions again and give remarks on potential extensions of the results to more general test functions and to the multidimensional case.
\label{rev:proof_extension_reasons_1}

\section{\texorpdfstring{On $\eps$}{Lg}-dependent limit theorems}  \label{sec:limit_theorems}

\subsection{An \texorpdfstring{$\eps$}{Lg}-dependent mean ergodic theorem}  \label{subsec:mean_ergodic_theorem}
Consider a probability space $\left( \Omega, \mathcal{F}, \Pb \right)$ with a filtration $\left( \mathcal{F}_t \right)_{t \in [0, \infty)}$ that satisfies the usual conditions, and a one-dimensional Brownian motion $W := \left( W(t), \mathcal{F}_t \right)_{t \in [0, \infty)}$ on that probability space. Assume further that there is another one-dimensional Brownian motion $\overline{W} := \left( \overline{W}(t), \mathcal{F}_t \right)_{t \in [0, \infty)}$ that is independent of $W$ and lives on the same probability space. This construction is feasible through a suitable extension of the original probability space.

Fix $x_0 \in \R$ once and for all. Consider the following one-dimensional stochastic differential equation (SDE) in $\R$, depending on a small parameter $\eps>0$,
\begin{equation} \label{eq:original_sde_eps}
    d X_\eps(t) = b_\eps(X_\eps(t)) dt + \sigma(X_\eps(t)) dW(t), \quad t > 0, \quad X_\eps(0) = x_0,
\end{equation}
with suitable Borel-measurable functions $b_\eps \colon \R \rightarrow \R$, $\sigma \colon \R \rightarrow \R$, and assume, for the moment, that there exists a unique, strong solution to this SDE on the given probability space. We suppose a setting where $X_\eps$ converges weakly to another stochastic process $X$ in $C([0,T];\R)$ as $\epstozero$. Here, $X$ is the unique, strong solution to the stochastic differential equation
\begin{equation} \label{eq:original_sde}
    d X(t) = b(X(t)) dt + \overline{\sigma}(X(t)) d\overline{W}(t), \quad t > 0, \quad X(0) = x_0,
\end{equation}
with Borel-measurable functions $b \colon \R \rightarrow \R$ and $\overline{\sigma} \colon \R \rightarrow \R$. Note that $\sigma \neq \overline{\sigma}$ in general. Under our Assumptions \ref{ass:assumptions_C} and \ref{ass:assumptions_MET}, which will be introduced soon, the following scale functions exist and are twice differentiable on $\R$ with strictly positive first derivative
\begin{align}       \label{eq:A-harmonic_functions}
    f_\eps(x) := \int_0^x \exp \left( -\int_0^z \frac{2 b_\eps(y)}{\sigma(y)^2} \, dy \right) \, dz, \quad
    f(x) := \int_0^x \exp \left( -\int_0^z \frac{2 b(y)}{\overline{\sigma}(y)^2} \, dy \right) \, dz.
\end{align}
These functions are harmonic with respect to the differential operators
\begin{align}       \label{eq:generators}
    \mathcal{A}_\eps := b_\eps \partial_x + \frac{1}{2} \sigma^2 \partial_{xx}, \quad
    \mathcal{A} := b \partial_x + \frac{1}{2} \overline{\sigma}^2 \partial_{xx},
\end{align}
respectively, that is $\mathcal{A}_\eps f_\eps = 0$ and $\mathcal{A} f = 0$. Setting $\xi_\eps := f_\eps \circ X_\eps$, $\xi := f \circ X$, and applying the Itô formula yields the transformed stochastic differential equations
\begin{align}
    \label{eq:transformed_sde_eps}
    &d \xi_\eps(t) = \frac{1}{\rho_\eps(\xi_\eps(t))} dW(t), \quad t > 0, \quad \xi_\eps(0) = f_\eps(x_0), \\[0.25cm]
    \label{eq:transformed_sde}
    &d \xi(t) = \frac{1}{\rho(\xi(t))} d\overline W(t), \quad t > 0, \quad \xi(0) = f(x_0).
\end{align}
Here, the functions $\rho_\eps$ and $\rho$ are given by
\begin{align}
    \rho_\eps(x) := \frac{1}{\sigma(g_\eps(x)) f'_\eps(g_\eps(x))}, \quad
    \rho(x) := \frac{1}{\overline{\sigma}(g(x)) f'(g(x))}, 
\end{align}
where $g$ and $g_\eps$ are the inverse functions to $f$ and $f_\eps$, respectively.

In the following we want to derive an "$\eps$-dependent" MET, i.e., under appropriate assumptions and for arbitrary initial conditions, we want to prove the validity of 
\begin{align}       \label{eq:eps_mean_ergodic_theorem}
    \E \left| \frac{1}{T_\eps} \int_0^{T_\eps} \varphi_\eps(\xi_\eps(t))  \, dt - \frac{1}{C_{\rho_\eps}} \int_\R \varphi_\eps(y) \rho_\eps(y)^2 \, dy \right|^2 &\rightarrow 0, \quad \epstozero,
\end{align}
where $\varphi_\eps \colon \R \to \R$ is a suitable measurable function, $C_{\rho_\eps} := \int_\R \rho_\eps(x)^2 \, dx$, $\eps>0$, and $T_\eps \rightarrow \infty$ as $\epstozero$. If we further assume the convergence
\begin{equation*}
    \frac{1}{C_{\rho_\eps}} \int_\R \varphi_\eps(y) \rho_\eps(y)^2 \, dy \rightarrow \frac{1}{C_{\rho}} \int_\R \varphi(y) \rho(y)^2 \, dy, \quad \epstozero,
\end{equation*}
for some function $\varphi \colon \R \to \R$ and $C_\rho := \int_\R \rho(x)^2 \, dx$, then it follows that
\begin{align}
    \E \left| \frac{1}{T_\eps} \int_0^{T_\eps} \varphi_\eps(\xi_\eps(t))  \, dt - \frac{1}{C_{\rho}} \int_\R \varphi(y) \rho(y)^2 \, dy \right|^2 &\rightarrow 0, \quad \epstozero.
\end{align}
\label{rev:proof_extension_reasons_2}
To prove \eqref{eq:eps_mean_ergodic_theorem} we will follow the material by Gikhman and Skorokhod, \cite{GS:1968}, on the one hand and by Sundar, \cite{S:1989}, on the other hand. They proved the MET for the case where no $\eps$-dependence is present and the limit is established for $T \rightarrow \infty$. We will carefully adapt and adequately combine the separate approaches whenever necessary.

In this work, we will constrict ourselves to the following assumptions:
\begin{assumptioncustom}{(C)}   \label{ass:assumptions_C}
    \hfill
    \begin{enumerate}[label=\roman*)]
    \item The functions $b$ and $\overline{\sigma}$ are locally Lipschitz-continuous on $\R$ and $\overline{\sigma}$ is strictly positive.
    \item The function $\sigma$ is locally Lipschitz-continuous on $\R$ and satisfies $a \leq \sigma \leq A$ on $\R$ for some $A, a>0$. Moreover, for each $\eps > 0$ and $N \in \N$ there exists a constant $L_N > 0$ such that for all $x, y \in [-N, N]$
    \begin{align*}
        |b_\eps(x) - b_\eps(y)| \leq \frac{L_N}{\eps^2} |x-y|, \quad
        |b_\eps(x)| \leq \frac{L_N}{\eps^2}(1 + |x|).
    \end{align*}
    \end{enumerate}
\end{assumptioncustom}
\begin{Remark}
Most homogenization results for SDEs are, to the best of our knowledge, established under global Lipschitz assumptions, cf. \cite{BLP:1978, P:1999, PV:2001, PS:2008, RX:2021}, but, according to \cite{PSV:1976} it is possible to have local assumptions as long as the solution to \eqref{eq:original_sde} has almost surely infinite explosion time.
\end{Remark}
\begin{assumptioncustom}{(MET)}   \label{ass:assumptions_MET}
    \hfill
    \begin{enumerate}[label=\roman*)]
    \item $\lim_{x \rightarrow \pm \infty} f(x) = \pm \infty$.
    \item $C_\rho < \infty.$
    \item There exist constants $c_f, c_\sigma, C_f, C_\sigma > 0$ such that for all $\eps > 0$ and on $\R$
    \begin{equation*}
        c_f f' \leq f'_\eps \leq C_f f', \qquad 
        c_\sigma \overline{\sigma} \leq \sigma \leq C_\sigma \overline{\sigma}.
    \end{equation*}
    \end{enumerate}
\end{assumptioncustom}
\noindent

We are not suggesting that these assumptions are minimal requirements for the following results to work. On the contrary, some of them, e.g., the pointwise bounds in \ref{ass:assumptions_MET} iii), can be rather restrictive and for certain proofs even superfluous. Indeed, the same results can be obtained if, for example, in \ref{ass:assumptions_MET} (iii) there exist two integrable functions $k_1, k_2$ such that $c_f k_1 \leq f'_\eps \leq C_f k_2$ on $\R$. But in order to keep the exposition clear without repeating similar arguments and altering assumptions over and over again, we will stick to them throughout this work. Note, however, that we do not impose global Lipschitz assumptions on the drift coefficients, as it was originally done in \cite{GS:1968} and \cite{S:1989}. We also need to emphasize the circumstance that the local Lipschitz and local linear growth constant in \ref{ass:assumptions_C} ii) explode as $\epstozero$ and we need to let the test functions in the integrals depend on $\eps$. This makes the analysis substantially more difficult.

\label{rev:proof_extension_reasons_3}
\begin{Remark}  \label{rem:existence_invariant_density}
    As it has been established in \cite[\S 18]{GS:1968} and \cite{S:1989}, but will also become apparent here, the invariant density of a stochastic process $X$ defined through an SDE such as \eqref{eq:original_sde} exists whenever its corresponding scale function $f(x) \rightarrow \pm \infty$ as $x \rightarrow \pm \infty$, and the function $\rho^2$ is integrable. In this case, the invariant density of $X$ is given by
    \begin{equation}
        \mu(x) := \frac{1}{Z \overline{\sigma}(x)^2} \exp \left( \int_0^x \frac{2 b(y)}{\overline{\sigma}(y)^2} \, dy \right), \quad Z := \int_\R \overline{\sigma}(x)^{-2} \exp \left( \int_0^x \frac{2 b(y)}{ \overline{\sigma}(y)^2} \, dy \right) \, dx.
    \end{equation}
    For fixed $\eps>0$ and under the same assumptions, a similar formula holds for the invariant density of $X_\eps$, namely
    \begin{equation} \label{eq:original_invariant_density_eps}
    \mu_\eps(x) := \frac{1}{Z_\eps \sigma(x)^2} \exp \left( \int_0^x \frac{2 b_\eps(y)}{\sigma(y)^2} \, dy \right), \quad Z_\eps := \int_\R \sigma(y)^{-2} \exp \left( \int_0^x \frac{2 b_\eps(y)}{\sigma(y)^2} \, dy \right) \, dx.
\end{equation}
\end{Remark}

Observe that an ordinary integral substitution gives
\begin{align*}
    \frac{1}{T_\eps} \int_0^{T_\eps} \varphi_\eps(X_\eps(t))  \, dt - \int_\R \varphi_\eps(x) \mu_\eps(x) \, dx
    = \frac{1}{T_\eps} \int_0^{T_\eps} \varphi_\eps(g_\eps(\xi_\eps(t)))  \, dt - \frac{1}{C_{\rho_\eps}} \int_\R \varphi_\eps(g_\eps(x)) \rho_\eps(x)^2 \, dx,
\end{align*}
indicating that it is indeed sufficient to prove the mean ergodic theorem for the solution of equation \eqref{eq:transformed_sde_eps}. We start with a lemma that summarizes some properties, which follow from Assumptions \ref{ass:assumptions_C}, \ref{ass:assumptions_MET}, and which will be used throughout what follows.

\begin{Lemma}   \label{lemma:A1_properties}
Assume \ref{ass:assumptions_C} and \ref{ass:assumptions_MET}.
\begin{enumerate}[label=\roman*)]
    \item For each of the SDEs defined in \eqref{eq:original_sde_eps}, \eqref{eq:original_sde}, \eqref{eq:transformed_sde_eps}, \eqref{eq:transformed_sde} with their given initial conditions there exists a unique, strong solution. Furthermore, for $\eps >0$ we have the moment inequality
    \begin{equation}   \label{eq:lyapunov_moment_inequality}
        \E \left| \xi_\eps(t) \right| \leq \left| f_\eps(x_0) \right| \exp(t), \quad t \geq 0.
    \end{equation}
    \item There exist constants $d_\rho, D_\rho>0$ such that $d_\rho \leq C_{\rho_\eps} \leq D_\rho $ for all $\eps > 0$.
    \item For every $\delta > 0$ there exists an $M > 0$ such that for all $\eps > 0$
    \begin{equation*}
        \int_{|y|>M} \rho_\eps(y)^2 \, dy < \delta.
    \end{equation*}
\end{enumerate}
\end{Lemma}
\begin{proof} \hfill \\
i) Strong uniqueness for \eqref{eq:transformed_sde_eps}, \eqref{eq:transformed_sde}, and inequality \eqref{eq:lyapunov_moment_inequality} are a consequence of \cite[Theorem 3.5]{KMN:2012}. Compare also with \cite[Example 3.10, Appendix A.2]{KMN:2012}. Strong uniqueness for \eqref{eq:original_sde_eps}, \eqref{eq:original_sde} follows from \cite[Proposition 5.5.13]{KS:1996}. Note that, according to Feller's test for explosions, the solution processes do not explode in finite time due to \ref{ass:assumptions_MET} i). \\[0.25cm]
ii) This follows by a simple integral substitution and \ref{ass:assumptions_MET} iii), namely
\begin{equation*}
    C_{\rho_\eps} = \int_\R \rho_\eps(y)^2 \, dy = \int_\R \frac{dy}{\sigma(y)^2 f_\eps'(y)} \leq \frac{1}{c_\sigma^2 c_f} \int_\R \rho(y)^2 \, dy =: D_\rho,
\end{equation*}
and, similarly, for the lower bound. \\[0.25cm]
iii) First note that by \ref{ass:assumptions_MET} iii) we have for any $\eps>0$
\begin{equation}    \label{eq:f_estimates}
    c_f f \leq f_\eps \leq C_f f, \quad \text{on } \R.
\end{equation}
Now let $\delta>0$ be arbitrary. Due to \ref{ass:assumptions_MET} ii) we can choose $M>0$ such that 
\begin{equation*}
    \int_{|y|>M/(c_f \vee C_f)} \rho(y)^2 \, dy < \delta c_\sigma^2 c_f,
\end{equation*}
where $c_f \vee C_f := \max\{c_f, C_f\}$. Using this, \eqref{eq:f_estimates}, and the strict monotonicity of $f, f_\eps, g,$ and $ g_\eps$, we can therefore estimate for all $\eps > 0$
\begin{align*}
    \int_{|y|>M} \rho_\eps(y)^2 \, dy &= \int_{|f_\eps(y)|>M} \frac{dy}{\sigma(y)^2 f'_\eps(y)} \leq \frac{1}{c_\sigma^2 c_f} \int_{|f(y)|>M/(c_f \vee C_f)} \frac{dy}{\overline{\sigma}(y)^2 f'(y)} \\[0.25cm]
    &= \frac{1}{c_\sigma^2 c_f} \int_{|y|>M/(c_f \vee C_f)} \rho(y)^2 \, dy < \delta. \tag*{\qedhere}
\end{align*}
\end{proof}

Under Assumptions \ref{ass:assumptions_C} and \ref{ass:assumptions_MET}, we define for any $x, y \in \R$ and $\eps > 0$ the stopping time $\tau^{y}\left(\xi_\eps^x\right) := \inf \{s>0 \, | \, \xi_\eps^x(s) = y\}$, i.e., the time that the process $\xi_\eps^x$ solving \eqref{eq:transformed_sde_eps} and starting in $\xi_\eps^x(0) = x$ requires to get to the point $y$ for the first time. Note that, by \cite[\S16, Theorem 1]{GS:1968},
\begin{equation*}
    \Pb\left(\sup_{t>0} \xi_\eps^x(t) = \infty\right) = \Pb\left(\inf_{t>0} \xi_\eps^x(t) = -\infty\right) = 1,
\end{equation*}
so that $\Pb(\tau^{y}\left(\xi_\eps^x\right) < \infty) = 1$, that is to say, $\xi_\eps^x$ is a recurrent process. Set also $\tau^{(a,b)} \left( \xi_\eps^x \right) := \inf \{s>0 \, | \, \xi_\eps^x(s) \notin (a,b) \}$ for any interval $(a,b) \subseteq \R$ and $x \in (a,b)$. By \cite[\S15, Corollary 15]{GS:1968} and some easy rearrangements of terms we get the following formula
\begin{equation*}
    \E \tau^{(a,b)} \left( \xi_\eps^x \right) = \frac{x-a}{b-a} \int_x^b  2(b-z)\rho_\eps(z)^2 \, dz + \frac{b-x}{b-a} \int_a^x  2(z-a)\rho_\eps(z)^2 \, dz.
\end{equation*}
The proof of the following result is found under \cite[\S18, Lemma 1]{GS:1968}. As $\eps>0$ is fixed in this case, nothing needs to be modified.
\begin{Lemma}   \label{lemma:random_time_lemma}
Under Assumptions \ref{ass:assumptions_C} and \ref{ass:assumptions_MET}, we have for any $x,y \in \R$ with $y < x$
\begin{equation}
    \E \tau^{y}\left(\xi_\eps^x\right) = 2(x-y) \int_x^\infty \rho_\eps(z)^2 \, dz + 2 \int_y^x  (z-y)\rho_\eps(z)^2 \, dz,
\end{equation}
and for $x < y$
\begin{equation}
    \E \tau^{y}\left(\xi_\eps^x\right) = 2(y-x) \int_{-\infty}^x \rho_\eps(z)^2 \, dz + 2\int_x^y  (y-z)\rho_\eps(z)^2 \, dz.
\end{equation}
\end{Lemma}

\begin{Proposition} \label{prop:semigroup_estimate}
Let $(\varphi_\eps)_{\eps>0}$ be a uniformly bounded sequence of Borel-measurable, real-valued functions on $\R$. Under Assumptions \ref{ass:assumptions_C} and \ref{ass:assumptions_MET}, it holds for any $x, y \in \R$ and $\eps > 0$
\begin{equation}    \label{eq:semigroup_estimate}
    \left|  \frac{1}{T_\eps} \int_0^{T_\eps} \E \varphi_\eps(\xi_\eps^x(t))  \, dt - \frac{1}{T_\eps} \int_0^{T_\eps} \E \varphi_\eps(\xi^y_\eps(t))  \, dt \right| \leq \frac{4 \sup_{\eps > 0} \| \varphi_\eps \|_{\infty} D_\rho |x-y|}{T_\eps}.
\end{equation}
\end{Proposition}
\begin{proof}
The proof is basically analogous to the original one, see \cite[\S18, Remark 1, Lemma 2]{GS:1968}, with the only difference being the estimate
\begin{equation*}
    \E \tau^{y}\left(\xi_\eps^x\right) \leq 2 D_\rho |x-y|,
\end{equation*}
which is an immediate consequence of Lemma \ref{lemma:random_time_lemma} and Lemma \ref{lemma:A1_properties} ii).
\end{proof}
Before we proceed, we have to introduce a sequence of nondecreasing stopping times that will aid us several times in the subsequent proofs. The introduction of this stopping time is mostly a technical inconvenience due to the local Lipschitz assumptions, but, nonetheless, required to make everything rigorous.\label{rev:martingale_notation} In the subsequent material, notation like $\left( M_\eps^z(t), \mathcal{F}_t\right)_{t \in [0, \infty)}$ stands for a real-valued stochastic process $\left( M_\eps^z(t) \right)_{t \in [0, \infty)}$ that is adapted to the filtration $\left( \mathcal{F}_t\right)_{t \in [0, \infty)}$. To carry out the localization procedure through the stopping times, fix $\eps>0$, $x \in \R$, and a bounded, Borel-measurable function $\varphi_\eps \colon \R \to \R$.  We define the following processes
\begin{align}
    K_\eps(t) &:= \int_0^t \varphi_\eps(\xi_\eps^{x}(s)) \, ds, \quad t \in [0, \infty), \\[0.25cm]
    M_\eps^z(t) &:= \int_0^t \sign(\xi_\eps^{x}(s) - z) \rho_\eps^{-1}(\xi_\eps^{x}(s)) \, dW(s), \quad t \in [0, \infty), \quad z \in \R,
\end{align}
where $\sign$ is defined to be left-continuous, i.e. $\sign = \mathds{1}_{(0, \infty)} - \mathds{1}_{(-\infty, 0]}$.
The appearing integrands, thus the integrals, are measurable and adapted to the filtration $\left( \mathcal{F}_t\right)_{t \in [0, \infty)}$. Observe that for any $t \geq 0$
\begin{equation*}
    \int_0^t \rho_\eps^{-2}(\xi_\eps^{x}(s)) \, ds \leq t \sup_{0 \leq s \leq t} \rho_\eps^{-2}(\xi_\eps^{x}(s)) < \infty, \quad \Pb\text{-a.s.},
\end{equation*}
by the continuity of $\rho_\eps^{-1}$ and $\xi_\eps^{x}$, so that $\left( M_\eps^z(t), \mathcal{F}_t\right)_{t \in [0, \infty)} \in \mathscr{M}^{\text{c, loc}}$ for each $z \in \R$ and $\eps > 0$, where $\mathscr{M}^{\text{c, loc}}$ is the space of continuous, local martingales starting in zero. Denote by $\langle M \rangle$ the quadratic variation process for $M \in \mathscr{M}^{\text{c, loc}}$. Given the boundedness of $\varphi_\eps$ and the fact that $d \langle M_\eps^z \rangle_s = \rho_\eps^{-2}(\xi_\eps^{x}(s)) \, ds$, a similar argument shows that
\begin{equation*}
    \left( I_\eps^z(t) := \int_0^t K_\eps(s) \, dM_\eps^z(s), \mathcal{F}_t\right)_{t \in [0, \infty)} \in \mathscr{M}^{\text{c, loc}},
\end{equation*}
for each $z \in \R$. Hence, there exists a nondecreasing sequence of stopping times $(\tau_n^{\varphi_\eps, x})_{n \in \N}$ of $\left( \mathcal{F}_t\right)_{t \in [0, \infty)}$ such that $\tau_n^{\varphi_\eps, x} \rightarrow \infty$ $\Pb$-a.s. as $n \rightarrow \infty$ and 
\begin{equation}   \label{eq:mean_ergodic_theorem_local_martingale}
    \left( M_\eps^z(t \wedge \tau_n^{\varphi_\eps, x}), \mathcal{F}_t\right)_{t \in [0, \infty)}, \quad \left( I_\eps^z(t \wedge \tau_n^{\varphi_\eps, x}), \mathcal{F}_t\right)_{t \in [0, \infty)} \in \mathscr{M}_2^{\text{c}},
\end{equation}
for each $n \in \N$, where $t \wedge \tau_n^{\varphi_\eps, x} := \min\{t, \tau_n^{\varphi_\eps, x}\}$ and $\mathscr{M}_2^{\text{c}}$ is the space of continuous, square-integrable martingales starting in zero. Note that the quadratic variation process of $I_\eps^z$ is given by 
\begin{equation*}
    \langle I_\eps^z \rangle_t = \int_0^t K_\eps(s)^2 \rho_\eps^{-2}(\xi_\eps^{x}(s)) \, ds, \quad t \geq 0.
\end{equation*}
In particular, it does not depend on $z \in \R$. Therefore, and because $\langle I_\eps^z \rangle$ and $\xi_\eps^{x}$ are continuous, we can additionally choose the localizing sequence $\tau_n^{\varphi_\eps, x}$ in such a way that 
\begin{equation}    \label{eq:localizing_bounds}
    \langle I_\eps^z \rangle_{\cdot \, \wedge \tau_n^{\varphi_\eps, x}} \leq n, \quad \text{and} \quad |\xi_\eps^{x}( \, \cdot \wedge \tau_n^{\varphi_\eps, x})| \leq n, \quad \Pb\text{-a.s.},
\end{equation}
for each $n \in \N$ and all $z \in \R$, where we use the notation $\xi_\eps^{x}( \, \cdot \wedge \tau_n^{\varphi_\eps, x})$ for the map $[0, \infty) \rightarrow \R; t \mapsto \xi_\eps^{x}( \, t \wedge \tau_n^{\varphi_\eps, x})$ elementwise on $\Omega$ and similarly for $\langle I_\eps^z \rangle_{\cdot \, \wedge \tau_n^{\varphi_\eps, x}}$.

In addition to that, recall the Meyer-Tanaka formula for a continuous, local martingale, say $Y(t) = Y(0) + M(t)$ with $(M(t), \mathcal F_t)_{t \in [0, \infty)} \in \mathscr{M}^{\text{c, loc}}$, that is, for $z \in \R$ and $t \geq 0$
\begin{equation}    \label{eq:meyer_tanaka_formula}
    |Y(t) - z| = |Y(0) - z| + \int_0^t \sign(Y(s) - z) \, d M(s) + \Lambda_t^z(Y), \quad \Pb \text{-a.s.}
\end{equation}
Here, $\Lambda_t^z(Y)$ is the local time process for $Y$, cf.\ \cite[Theorem 3.7.1]{KS:1996}. Recall also the occupation density formula, which gives us for every Borel-measurable function $k \colon \R \rightarrow [0, \infty)$ the identity
\begin{equation}    \label{eq:occupation_density_formula}
    \int_0^t k(Y_s) \, d\langle M \rangle_s = \int_\R k(z) \Lambda_t^z(Y) \, dz, \quad t \geq 0, \quad \Pb \text{-a.s.}
\end{equation}

We will now utilize these formulas and the introduced sequence of stopping times to prove a result similar to Proposition \ref{prop:semigroup_estimate}. The proof idea is borrowed from \cite{S:1989}, but note that we do not assume global Lipschitz assumptions here and we need an explicit estimate rather than just a qualitative convergence statement.

\begin{Proposition}     \label{prop:semigroup_estimate_mod}
Let $(\varphi_\eps)_{\eps>0}$ be a uniformly bounded sequence of Borel-measurable, nonnegative functions on $\R$ with $\cup_{\eps > 0} \, \mathrm{supp}(\varphi_\eps) \subset [-R, R]$ for some $R>0$. Under Assumptions \ref{ass:assumptions_C} and \ref{ass:assumptions_MET}, it holds for any $x, y \in \R$ and $\eps > 0$
\begin{equation}    \label{eq:semigroup_estimate_t}
    \left| \frac{2}{T_\eps^2} \int_0^{T_\eps} t \, \E \varphi_\eps(\xi_\eps^x(t)) \, dt - \frac{2}{T_\eps^2} \int_0^{T_\eps} t \, \E \varphi_\eps(\xi^y_\eps(t)) \, dt \right| \leq \frac{8 \sup_{\eps > 0} \| \varphi_\eps \|_{\infty} D_\rho |x-y|}{T_\eps}.
\end{equation}
\end{Proposition}
\begin{proof}
Assume without loss of generality that $x \geq y$. Using Fubini's theorem, we first have
\begin{align}   \label{eq:start_proof_prop_25}
    \begin{aligned} 
    &\frac{2}{T_\eps^2} \int_0^{T_\eps} t \, \E \varphi_\eps(\xi_\eps^x(t)) \, dt - \frac{2}{T_\eps^2} \int_0^{T_\eps} t \, \E \varphi_\eps(\xi^y_\eps(t)) \, dt \\[0.25cm]
    = \, &\frac{2}{T_\eps^2} \int_0^{T_\eps} \int_s^{T_\eps} \E \varphi_\eps(\xi_\eps^x(t)) - \E \varphi_\eps(\xi^y_\eps(t)) \, dt \, ds \\[0.25cm]
    = \, &\frac{2}{T_\eps} \int_0^{T_\eps} \E \varphi_\eps(\xi_\eps^x(t)) - \E \varphi_\eps(\xi^y_\eps(t)) \, dt  - \frac{2}{T_\eps^2} \int_0^{T_\eps} \E \int_0^s \varphi_\eps(\xi_\eps^x(t)) - \varphi_\eps(\xi^y_\eps(t)) \, dt \, ds.    
    \end{aligned}
\end{align}
Define the sequence of stopping times $\tau_n := \tau_n^{\varphi_\eps, x} \wedge \tau_n^{\varphi_\eps, y}$ through the previously constructed stopping times. With the occupation density formula \eqref{eq:occupation_density_formula} we can write for any $s \in [0, \infty)$
\begin{align}  \label{eq:exp_formula_prop_25}
    \E \int_0^{s \wedge \tau_n} \varphi_\eps(\xi_\eps^x(t)) - \varphi_\eps(\xi^y_\eps(t)) \, dt
    = \E \int_\R \varphi_\eps(z) \rho_\eps(z)^2 \left[ L_{s \wedge \tau_n}^z(\xi_\eps^{x}) - L_{s \wedge \tau_n}^z(\xi_\eps^{y}) \right] \, dz.
\end{align}
The Meyer-Tanaka formula \eqref{eq:meyer_tanaka_formula} gives $\Pb$-a.s.
\begin{equation*}
    L_{s \wedge \tau_n}^z(\xi_\eps^{x}) = |\xi_\eps^{x}(s \wedge \tau_n) - z| - |x - z| - \int_0^{s \wedge \tau_n} \sign(\xi_\eps^{x}(r) - z) \rho_\eps^{-1}(\xi_\eps^{x}(r)) \, dW(r),
\end{equation*}
and, similarly, for $L_{s \wedge \tau_n}^z(\xi_\eps^{y})$. The respective stochastic integrals in these formulas are martingales, cf.\ \eqref{eq:mean_ergodic_theorem_local_martingale}, starting in zero and, thus, have expectation zero. Hence
\begin{equation}
    \E \left[ L_{s \wedge \tau_n}^z(\xi_\eps^{x}) - L_{s \wedge \tau_n}^z(\xi_\eps^{y}) \right] = \E|\xi_\eps^{x}(s \wedge \tau_n) - z| - |x - z| - \E|\xi_\eps^{y}(s \wedge \tau_n) - z| + |y - z|.
\end{equation}
With the uniqueness of the solution of the SDE \eqref{eq:transformed_sde_eps} and $x \geq y$, a stopping time argument shows that $\xi_\eps^{x}(r) \geq \xi_\eps^{y}(r)$ for all $r \in [0, \infty)$, $\Pb$-a.s. Using this fact, a simple case distinction shows that
\begin{align*}
    \E|\xi_\eps^{x}(s \wedge \tau_n) - z| - \E|\xi_\eps^{y}(s \wedge \tau_n) - z|
    &\leq \E \left[\xi_\eps^{x}(s \wedge \tau_n) - \xi_\eps^{y}(s \wedge \tau_n)\right] = x - y, \\[0.25cm]
    |y - z| - |x - z| &\leq x-y,
\end{align*}
where we used the martingale property in the first line. Therefore,
\begin{equation}   \label{eq:estimate_local_time_difference}
    \E \left[ L_{s \wedge \tau_n}^z(\xi_\eps^{x}) - L_{s \wedge \tau_n}^z(\xi_\eps^{y}) \right] \leq 2(x-y).
\end{equation}
Observe that by \eqref{eq:localizing_bounds}, the compact support of $\varphi_\eps$, and the Burkholder-Davis-Gundy inequalities we can upper bound
\begin{align*}
    &\E \int_\R |\varphi_\eps(z)| \rho_\eps(z)^2 L_{s \wedge \tau_n}^z(\xi_\eps^{x}) \, dz \\[0.25cm]
    \leq \, & \E \int_\R |\varphi_\eps(z)| \rho_\eps(z)^2 \left[ |\xi_\eps^{x}(s \wedge \tau_n) - z| + |x - z| + \left| \int_0^{s \wedge \tau_n} \sign(\xi_\eps^{x}(r) - z) \rho_\eps^{-1}(\xi_\eps^{x}(r)) \, dW(r) \right| \right] \, dz \\[0.25cm]
    \leq \, & \| \varphi_\eps \|_\infty \int_{[-R, R]} \rho_\eps(z)^2 \left[ \E |\xi_\eps^{x}(s \wedge \tau_n) | + |x| + 2|z| + \E \left| \int_0^{s \wedge \tau_n} \sign(\xi_\eps^{x}(r) - z) \rho_\eps^{-1}(\xi_\eps^{x}(r)) \, dW(r) \right| \right] \, dz \\[0.25cm]
    \leq \, & \| \varphi_\eps \|_\infty D_\rho \left[ n + |x| + 2R + B \E \langle I_\eps^z \rangle_{s \wedge \tau_n}^{1/2} \right]
    \leq \| \varphi_\eps \|_\infty D_\rho \left[ n + |x| + 2R + B n^{1/2} \right] < \infty,
\end{align*}
where $B>0$ is the universal constant coming from the Burkholder-Davis-Gundy inequalities. A completely analog estimate is obtained when $L_{s \wedge \tau_n}^z(\xi_\eps^{x})$ is replaced with $L_{s \wedge \tau_n}^z(\xi_\eps^{y})$. Hence, we can use Fubini's theorem together with \eqref{eq:estimate_local_time_difference} to obtain the estimate
\begin{equation*}
    \E \int_\R \varphi_\eps(z) \rho_\eps(z)^2 \left[ L_{s \wedge \tau_n}^z(\xi_\eps^{x}) - L_{s \wedge \tau_n}^z(\xi_\eps^{y}) \right] \, dz \leq 2(x-y) \int_\R \varphi_\eps(z) \rho_\eps(z)^2 \, dz.
\end{equation*}
Letting $n \rightarrow \infty$, the dominated convergence theorem implies in \eqref{eq:exp_formula_prop_25} for any $s \in [0, \infty)$ 
\begin{equation*}
    \E \int_0^s \varphi_\eps(\xi_\eps^x(t)) - \varphi_\eps(\xi^y_\eps(t)) \, dt \leq 2(x-y) \int_\R \varphi_\eps(z) \rho_\eps(z)^2 \, dz \leq 2(x-y) \sup_{\eps > 0} \| \varphi_\eps \|_{\infty} D_\rho.
\end{equation*}
Therefore, we conclude from \eqref{eq:start_proof_prop_25} and Proposition \ref{prop:semigroup_estimate}
\begin{equation*}
    \left| \frac{2}{T_\eps^2} \int_0^{T_\eps} \E \varphi_\eps(\xi_\eps^x(t)) t \, dt - \frac{2}{T_\eps^2} \int_0^{T_\eps} \E \varphi_\eps(\xi^y_\eps(t)) t \, dt \right|
    \leq \frac{8 \sup_{\eps > 0} \| \varphi_\eps \|_{\infty} D_\rho (x-y)}{T_\eps}.      \tag*{\qedhere}
\end{equation*}
\end{proof}
The proof of the next stationarity result can be obtained, once again, by the same arguments as in the original material in \cite[\S18]{GS:1968}.
\begin{Corollary}
Under Assumptions \ref{ass:assumptions_C} and \ref{ass:assumptions_MET}, we have for any Borel-measurable function $\varphi_\eps \colon \R \to \R$ with $\int_\R |\varphi_\eps(x)| \rho_\eps(x)^2 \, dx < \infty$ and for $t\geq0$
\begin{equation}    \label{eq:stationarity}
    \int_\R \E \varphi_\eps(\xi^y_\eps(t)) \rho_\eps(y)^2 \, dy = \int_\R \varphi_\eps(y) \rho_\eps(y)^2 \, dy.
\end{equation}
\end{Corollary}

The following lemma gives us another key ingredient for our sought-after main result.
\begin{Lemma}   \label{lemma:semigroup_ergodicity}
Let $(\varphi_\eps)_{\eps>0}$ be a uniformly bounded sequence of Borel-measurable, nonnegative functions on $\R$. Under Assumptions \ref{ass:assumptions_C} and \ref{ass:assumptions_MET}, it holds for $x_\eps := f_\eps(x_0)$
\begin{align} \label{eq:semigroup_ergodicity}
\begin{aligned}
    &\left|  \frac{1}{T_\eps} \int_0^{T_\eps} \E \varphi_\eps(\xi_\eps^{x_{\eps}}(t))  \, dt - \frac{1}{C_{\rho_\eps}} \int_\R \varphi_\eps(y) \rho_\eps(y)^2 \, dy \right| \rightarrow 0, \quad \epstozero,
\end{aligned}
\end{align}
and, similarly, if $\cup_{\eps > 0} \mathrm{supp}(\varphi_\eps) \subset [-R, R]$ for some $R>0$, then
\begin{align}
\begin{aligned}
    &\left|  \frac{2}{T_\eps^2} \int_0^{T_\eps} t \, \E \varphi_\eps(\xi_\eps^{x_{\eps}}(t)) \, dt - \frac{1}{C_{\rho_\eps}} \int_\R \varphi_\eps(y) \rho_\eps(y)^2 \, dy \right| \rightarrow 0, \quad \epstozero.
\end{aligned}
\end{align}
\end{Lemma}
\begin{proof}
We will only prove the first claim \eqref{eq:semigroup_ergodicity} as the proof of the second claim works the same. For that second claim we would only have to replace estimate \eqref{eq:semigroup_estimate} with \eqref{eq:semigroup_estimate_t} in the following. 

To prove claim \eqref{eq:semigroup_ergodicity}, we first obtain from equation \eqref{eq:stationarity}
\begin{align*}
    &\left|  \frac{1}{T_\eps} \int_0^{T_\eps} \E \varphi_\eps(\xi_\eps^{x_{\eps}}(t))  \, dt - \frac{1}{C_{\rho_\eps}} \int_\R \varphi_\eps(y) \rho_\eps(y)^2 \, dy \right| \\[0.25cm]
    = &\frac{1}{C_{\rho_\eps}} \left|  \frac{1}{T_\eps} \int_0^{T_\eps} \E \varphi_\eps(\xi_\eps^{x_{\eps}}(t))  \, dt \; \int_\R \rho_\eps(y)^2 \, dy - \frac{1}{T_\eps} \int_0^{T_\eps} \int_\R \E \varphi_\eps(\xi^y_\eps(t)) \rho_\eps(y)^2 \, dy \, dt \right| \\[0.25cm]
    = &\frac{1}{C_{\rho_\eps}} \left|  \int_\R \left( \frac{1}{T_\eps} \int_0^{T_\eps} \E \varphi_\eps(\xi_\eps^{x_{\eps}}(t)) - \E \varphi_\eps(\xi^y_\eps(t)) \, dt \right) \rho_\eps(y)^2 \, dy \right| \\[0.25cm]
    \leq &\frac{1}{d_\rho} \left[ 2\sup_{\eps > 0} \| \varphi_\eps \|_{\infty} \int_{|x_{\eps}-y|>T_{\eps}} \rho_\eps(y)^2 \, dy + \int_{|x_{\eps}-y| \leq T_{\eps}} \left| \frac{1}{T_\eps} \int_0^{T_\eps} \E \varphi_\eps(\xi_\eps^{x_{\eps}}(t)) - \E \varphi_\eps(\xi^y_\eps(t)) \, dt \right| \rho_\eps(y)^2 \, dy \right].
\end{align*}
For the second term inside the brackets we use estimate \eqref{eq:semigroup_estimate} to get
\begin{align}
\begin{aligned}
    &\left|  \frac{1}{T_\eps} \int_0^{T_\eps} \E \varphi_\eps(\xi_\eps^{x_{\eps}}(t))  \, dt - \frac{1}{C_{\rho_\eps}} \int_\R \varphi_\eps(y) \rho_\eps(y)^2 \, dy \right| \\[0.25cm]
    \leq &\frac{2 \sup_{\eps > 0} \| \varphi_\eps \|_{\infty} }{d_\rho} \left[ \int_{|x_{\eps}-y|>T_{\eps}} \rho_\eps(y)^2 \, dy + \frac{2D_\rho}{T_\eps} \int_{|x_{\eps}-y| \leq T_{\eps}} |x_{\eps}-y| \rho_\eps(y)^2 \, dy \right].
\end{aligned}
\end{align}
Now let $\delta>0$ be arbitrary. By Lemma \ref{lemma:A1_properties} iii) we can choose $M>0$ such that for all $\eps>0$
\begin{equation*}
    \int_{|y|>M} \rho_\eps(y)^2 \, dy < \delta.
\end{equation*}
Therefore, for sufficiently small $\eps > 0$ we have with \eqref{eq:f_estimates}
\begin{align*}
    \int_{|x_\eps-y|>T_\eps} \rho_\eps(y)^2 \, dy \leq \int_{|y|>T_\eps-|x_\eps|} \rho_\eps(y)^2 \, dy \leq \int_{|y|>T_\eps - (c_f \vee C_f) |f(x_0)|}\rho_\eps(y)^2 \, dy < \delta,
\end{align*}
and 
\begin{align*}
    &\frac{1}{T_\eps} \int_{|x_\eps-y| \leq T_\eps} |x_\eps-y| \rho_\eps(y)^2 \, dy \\[0.25cm]
    = &\frac{1}{T_\eps} \left[ \int_{|x_\eps-y| \leq T_\eps, \; |y| \leq M} |x_\eps-y| \rho_\eps(y)^2 \, dy + \int_{|x_\eps-y| \leq T_\eps, \; |y| > M} |x_\eps-y| \rho_\eps(y)^2 \, dy \right] \\[0.25cm]
    \leq &\frac{\left[(c_f \vee C_f) |f(x_0)| + M\right]D_\rho}{T_\eps} + \int_{|y| > M} \rho_\eps(y)^2 \, dy \\[0.25cm]
    \leq &\frac{\left[(c_f \vee C_f) |f(x_0)| + M\right]D_\rho}{T_\eps} + \delta.
\end{align*}
Thus
\begin{equation*}
    \limsup_{\epstozero} \left|  \frac{1}{T_\eps} \int_0^{T_\eps} \E \varphi_\eps(\xi_\eps^{x_{\eps}}(t))  \, dt - \frac{1}{C_{\rho_\eps}} \int_\R \varphi_\eps(y) \rho_\eps(y)^2 \, dy \right| \leq \frac{2 \sup_{\eps > 0} \| \varphi_\eps \|_{\infty} (1 + 2D_\rho)}{d_\rho} \delta.
\end{equation*}
Letting $\delta \rightarrow 0$ yields the claim.
\end{proof}

The next lemma will be used twice in the proof of the MET. Compared to the results in \cite{S:1989}, we again resort to localization arguments with the previously introduced sequence of stopping times.
\begin{Lemma}   \label{lemma:double_integral_claim}
Let $\varphi_\eps, \psi_\eps$ be bounded, Borel-measurable, nonnegative functions such that $\mathrm{supp}(\varphi_\eps) \subset [-R_1, R_1]$ and $\mathrm{supp}(\psi_\eps) \subset [-R_2, R_2]$ for some $R_1, R_2>0$. Set $\tau_n := \tau_n^{\varphi_\eps, x_\eps}$. Under Assumptions \ref{ass:assumptions_C} and \ref{ass:assumptions_MET}, we have for all $\eps > 0$ and $n \in \N$
\begin{align} \label{eq:double_integral_claim}
\begin{aligned}
    &\E \int_0^{T_\eps \wedge \tau_n} \varphi_\eps(\xi_\eps^{x_{\eps}}(t)) \int_{t}^{T_\eps \wedge \tau_n} \, \psi_\eps(\xi_\eps^{x_{\eps}}(s)) \, ds \, dt \\[0.25cm]
    = \, &\E \int_\R \psi_\eps(z) \rho_\eps(z)^2 \int_0^{T_\eps \wedge \tau_n} \varphi_\eps(\xi_\eps^{x_{\eps}}(t) \left[ |\xi_\eps^{x_\eps}(T_\eps \wedge \tau_n) - z| - |\xi_\eps^{x_\eps}(t) - z| \right] dt \, dz.
\end{aligned}    
\end{align}
\end{Lemma}
\begin{proof}
Applying the occupation density formula \eqref{eq:occupation_density_formula} gives $\Pb$-a.s.
\begin{align*}
    &\int_0^{T_\eps \wedge \tau_n} \varphi_\eps(\xi_\eps^{x_{\eps}}(t)) \int_{t}^{T_\eps \wedge \tau_n} \, \psi_\eps(\xi_\eps^{x_{\eps}}(s)) \, ds \, dt \\[0.25cm]
    = &\int_0^{T_\eps \wedge \tau_n} \varphi_\eps(\xi_\eps^{x_{\eps}}(t)) \int_\R \psi_\eps(z) \rho_\eps(z)^2 \left[ L_{T_\eps \wedge \tau_n}^z(\xi_\eps^{x_\eps}) -  L_{t}^z(\xi_\eps^{x_\eps}) \right] dz \, dt \\[0.25cm]
    = &\int_\R \psi_\eps(z) \rho_\eps(z)^2 \int_0^{T_\eps \wedge \tau_n} \varphi_\eps(\xi_\eps^{x_{\eps}}(t)) \left[ L_{T_\eps \wedge \tau_n}^z(\xi_\eps^{x_\eps}) -  L_{t}^z(\xi_\eps^{x_\eps}) \right] dt \, dz,
\end{align*}
where we used Tonelli's theorem in the last equality since $L^z_t(\xi_\eps^{x_\eps})$ is nondecreasing in $t$. The Meyer-Tanaka formula \eqref{eq:meyer_tanaka_formula} yields $\Pb$-a.s.
\begin{align*}
    L_{T_\eps \wedge \tau_n}^z(\xi_\eps^{x_\eps}) -  L_{t}^z(\xi_\eps^{x_\eps})
    &= |\xi_\eps^{x_\eps}(T_\eps \wedge \tau_n) - z| - |\xi_\eps^{x_\eps}(t) - z| - \int_{t}^{T_\eps \wedge \tau_n} \sign(\xi_\eps^{x_\eps}(s) - z) \rho_\eps^{-1}(\xi_\eps^{x_\eps}(s)) \, dW(s) \\[0.25cm]
    &= |\xi_\eps^{x_\eps}(T_\eps \wedge \tau_n) - z| - |\xi_\eps^{x_\eps}(t) - z| - \left[ M_\eps^z(T_\eps \wedge \tau_n) - M_\eps^z(t) \right],
\end{align*}
so that $\Pb$-a.s.
\begin{align*}
    &\int_0^{T_\eps \wedge \tau_n} \varphi_\eps(\xi_\eps^{x_{\eps}}(t)) \int_{t}^{T_\eps \wedge \tau_n} \, \psi_\eps(\xi_\eps^{x_{\eps}}(s)) \, ds \, dt \\[0.25cm]
    = \, &\int_\R \psi_\eps(z) \rho_\eps(z)^2 \Bigg[ \int_0^{T_\eps \wedge \tau_n} \varphi_\eps(\xi_\eps^{x_{\eps}}(t)) \left[ |\xi_\eps^{x_\eps}(T_\eps \wedge \tau_n) - z| - |\xi_\eps^{x_\eps}(t) - z| \right] dt \\[0.25cm]
    &-  \int_0^{T_\eps \wedge \tau_n} \varphi_\eps(\xi_\eps^{x_{\eps}}(t)) \left[ M_\eps^z(T_\eps \wedge \tau_n) - M_\eps^z(t) \right] dt \Bigg] dz,
\end{align*}
since $M_\eps^z$, $\xi_\eps^{x_\eps}$ are continuous processes $\Pb$-a.s. and $\varphi_\eps$ is bounded with compact support. Stochastic integration by parts implies $\Pb$-a.s.
\begin{align*}
    \int_0^{T_\eps \wedge \tau_n} K_\eps(t) \, dM_\eps^z(t)
    &= K_\eps(T_\eps \wedge \tau_n) M_\eps^z(T_\eps \wedge \tau_n) - \int_0^{T_\eps \wedge \tau_n} M_\eps^z(t) \varphi_\eps(\xi_\eps^{x_{\eps}}(t)) \, dt \\[0.25cm]
    &= \int_0^{T_\eps \wedge \tau_n} \varphi_\eps(\xi_\eps^{x_{\eps}}(t)) \left[ M_\eps^z(T_\eps \wedge \tau_n) - M_\eps^z(t) \right] \, dt,
\end{align*}
that is, we have $\Pb$-a.s.
\begin{align*}
    &\int_0^{T_\eps \wedge \tau_n} \varphi_\eps(\xi_\eps^{x_{\eps}}(t)) \int_{t}^{T_\eps \wedge \tau_n} \, \psi_\eps(\xi_\eps^{x_{\eps}}(s)) \, ds \, dt \\[0.25cm]
    = \, &\int_\R \psi_\eps(z) \rho_\eps(z)^2 \Bigg[ \int_0^{T_\eps \wedge \tau_n} \varphi_\eps(\xi_\eps^{x_{\eps}}(t)) \left[ |\xi_\eps^{x_\eps}(T_\eps \wedge \tau_n) - z| - |\xi_\eps^{x_\eps}(t) - z| \right] dt \\[0.25cm]
    &- \int_0^{T_\eps \wedge \tau_n} K_\eps(t) \, dM_\eps^z(t) \Bigg] dz.
\end{align*}
We require integrability to apply linearity and Fubini's theorem with respect to $\E$ in what follows. For this, observe that by \eqref{eq:lyapunov_moment_inequality} and \eqref{eq:localizing_bounds}
\begin{align*}
    &\E \int_0^{T_\eps \wedge \tau_n} \varphi_\eps(\xi_\eps^{x_{\eps}}(t)) \int_\R \psi_\eps(z) \rho_\eps(z)^2 \left[ |\xi_\eps^{x_\eps}(T_\eps \wedge \tau_n) - z| + |\xi_\eps^{x_\eps}(t) - z| \right] dz \, dt \\[0.25cm]
    \leq \, & T_\eps \| \varphi_\eps \|_\infty \| \psi_\eps \|_\infty D_\rho \left( \E | \xi_\eps^{x_\eps}(T_\eps \wedge \tau_n) | + 2 R_2 + \frac{1}{T_\eps}\int_0^{T_\eps} \E |\xi_\eps^{x_\eps}(t)| \, dt \right) \\[0.25cm]
    \leq \, & T_\eps \| \varphi_\eps \|_\infty \| \psi_\eps \|_\infty D_\rho \left( n + 2 R_2 + \frac{|f_\eps(x_0)|}{T_\eps}\int_0^{T_\eps} \exp(t) \, dt \right) < \infty,
\end{align*}
and using the Burkholder-Davis-Gundy inequalities we can upper bound with a universal constant $B>0$ and \eqref{eq:localizing_bounds}
\begin{align*}
    \E \int_\R \psi_\eps(z) \rho_\eps(z)^2 \left| \int_0^{T_\eps \wedge \tau_n} K_\eps(t) \, dM_\eps^z(t) \right| \, dz 
    &\leq B \int_\R \psi_\eps(z) \rho_\eps(z)^2 \, \E \langle I_\eps^z \rangle_{T_\eps \wedge \tau_n}^{1/2} \, dz \\[0.25cm]
    &\leq B \| \psi_\eps \|_\infty D_\rho n^{1/2} < \infty.
\end{align*}
Hence,
\begin{align*}
    &\E \int_0^{T_\eps \wedge \tau_n} \varphi_\eps(\xi_\eps^{x_{\eps}}(t)) \int_{t}^{T_\eps \wedge \tau_n} \, \psi_\eps(\xi_\eps^{x_{\eps}}(s)) \, ds \, dt \\[0.25cm]
    = \, & \int_\R \psi_\eps(z) \rho_\eps(z)^2 \Bigg[ \E \int_0^{T_\eps \wedge \tau_n} \varphi_\eps(\xi_\eps^{x_{\eps}}(t)) \left[ |\xi_\eps^{x_\eps}(T_\eps \wedge \tau_n) - z| - |\xi_\eps^{x_\eps}(t) - z| \right] dt \\[0.25cm]
    &- \E \int_0^{T_\eps \wedge \tau_n} K_\eps(t) \, dM_\eps^z(t) \Bigg] dz.
\end{align*}
Finally, we use the fact that the last displayed expectation is zero due to \eqref{eq:mean_ergodic_theorem_local_martingale} to reach the claim.
\end{proof}
Using these preparatory results, we are now in a position to state and prove the MET.
\begin{Theorem}     \label{theo:eps_mean_ergodic_theorem}
Let $(\varphi_\eps)_{\eps>0}$ be a uniformly bounded sequence of Borel-measurable, real-valued functions on $\R$. Under Assumptions \ref{ass:assumptions_C} and \ref{ass:assumptions_MET}, it holds for $x_\eps = f_\eps(x_0)$ and $T_\eps \rightarrow \infty$ as $\epstozero$
\begin{equation} \label{eq:eps_mean_ergodic_theorem_transf}
    \E \left| \frac{1}{T_\eps} \int_0^{T_\eps} \varphi_\eps(\xi_\eps^{x_{\eps}}(t))  \, dt - \frac{1}{C_{\rho_\eps}} \int_\R \varphi_\eps(y) \rho_\eps(y)^2 \, dy \right|^2 \rightarrow 0, \quad \epstozero.
\end{equation}
\end{Theorem}
\begin{Remark}  \label{rem:mean_ergodic_theorem_original_quantities}
In terms of the original quantities, i.e., with the process $X_\eps^{x_0}$, that solves \eqref{eq:original_sde_eps} and starts in $x_0$, and the invariant density $\mu_\eps$ in \eqref{eq:original_invariant_density_eps}, the limit result \eqref{eq:eps_mean_ergodic_theorem_transf} reads as
\begin{equation} \label{eq:eps_mean_ergodic_theorem_orig}
    \E \left| \frac{1}{T_\eps} \int_0^{T_\eps} \varphi_\eps(X_\eps^{x_0}(t))  \, dt - \int_\R \varphi_\eps(y) \mu_\eps(y) \, dy \right|^2 \rightarrow 0, \quad \epstozero.
\end{equation}
\end{Remark}
\begin{proof}
We will first prove the result for the case where the functions $\varphi_\eps$ are nonnegative with $\cup_{\eps > 0} \; \mathrm{supp}(\varphi_\eps) \subset [-R, R]$ for some $R>0$. Note that we can write
\begin{align*}
    &\E \left| \frac{1}{T_\eps} \int_0^{T_\eps} \varphi_\eps(\xi_\eps^{x_{\eps}}(t))  \, dt - \frac{1}{C_{\rho_\eps}} \int_\R \varphi_\eps(y) \rho_\eps(y)^2 \, dy \right|^2 \\[0.25cm]
    =\; &\E \left[ \left| \frac{1}{T_\eps} \int_0^{T_\eps} \varphi_\eps(\xi_\eps^{x_{\eps}}(t))  \, dt \right|^2 \right] - \left( \frac{1}{C_{\rho_\eps}} \int_\R \varphi_\eps(y) \rho_\eps(y)^2 \, dy \right)^2 \\[0.25cm]
    &- \E \left( \frac{1}{T_\eps} \int_0^{T_\eps} \varphi_\eps(\xi_\eps^{x_{\eps}}(t))  \, dt - \frac{1}{C_{\rho_\eps}} \int_\R \varphi_\eps(y) \rho_\eps(y)^2 \, dy \right) \frac{2}{C_{\rho_\eps}} \int_\R \varphi_\eps(y) \rho_\eps(y)^2 \, dy.
\end{align*}
It is enough to prove 
\begin{equation}    \label{eq:mean_ergodic_theorem_limsup_claim}
    \limsup_{\epstozero} \E \left[ \left| \frac{1}{T_\eps} \int_0^{T_\eps} \varphi_\eps(\xi_\eps^{x_{\eps}}(t))  \, dt \right|^2 \right] - \left( \frac{1}{C_{\rho_\eps}} \int_\R \varphi_\eps(y) \rho_\eps(y)^2 \, dy \right)^2 \leq 0,
\end{equation}
because then \eqref{eq:eps_mean_ergodic_theorem_transf} follows from the assumptions and Lemma \ref{lemma:semigroup_ergodicity}. For this, let $\delta>0$ be arbitrary. Choose $M>0$ such that for all $\eps>0$
\begin{equation*}
    \int_{|y|>M} \rho_\eps(y)^2 \, dy < \delta.
\end{equation*}
With the sequence of stopping times $(\tau_n)_{n \in \N}$ from the preceding proof of Lemma \ref{lemma:double_integral_claim} we have for all $\eps > 0$ and $n \in \N$
\begin{align}   \label{eq:double_integral_claim_1}
\begin{aligned}
    &\E \int_0^{T_\eps \wedge \tau_n} \varphi_\eps(\xi_\eps^{x_{\eps}}(t)) \int_{t}^{T_\eps \wedge \tau_n} \, \varphi_\eps(\xi_\eps^{x_{\eps}}(s)) \, ds \, dt \\[0.25cm]
    = \, &\E \int_\R \varphi_\eps(z) \rho_\eps(z)^2 \int_0^{T_\eps \wedge \tau_n} \varphi_\eps(\xi_\eps^{x_{\eps}}(t)) \left[ |\xi_\eps^{x_\eps}(T_\eps \wedge \tau_n) - z| - |\xi_\eps^{x_\eps}(t) - z| \right] dt \, dz,
\end{aligned}
\end{align}
and
\begin{align}   \label{eq:double_integral_claim_2}
\begin{aligned}
    &\E \int_0^{T_\eps \wedge \tau_n} \varphi_\eps(\xi_\eps^{x_{\eps}}(t)) \int_{t}^{T_\eps \wedge \tau_n} \,  \mathds{1}_{[-M, M]}(\xi_\eps^{x_{\eps}}(s)) \, ds \, dt \\[0.25cm]
    = \, &\E \int_{|y| \leq M} \rho_\eps(y)^2 \int_0^{T_\eps \wedge \tau_n} \varphi_\eps(\xi_\eps^{x_{\eps}}(t)) \left[ |\xi_\eps^{x_\eps}(T_\eps \wedge \tau_n) - y| - |\xi_\eps^{x_\eps}(t) - y| \right] dt \, dy.
\end{aligned}
\end{align}
First observe that
\begin{align*}
    &\int_\R \rho_\eps(y)^2 \, dy \, \E \left[ \left| \int_0^{T_\eps \wedge \tau_n} \varphi_\eps(\xi_\eps^{x_{\eps}}(t))  \, dt \right|^2 \right] \\[0.25cm]
    \leq \, &\left(\sup_{\eps>0} \| \varphi_\eps \|_\infty\right)^2 \E(T_\eps \wedge \tau_n)^2 \delta + \int_{|y| \leq M} \rho_\eps(y)^2 \, dy \, \E \left[ \left| \int_0^{T_\eps \wedge \tau_n} \varphi_\eps(\xi_\eps^{x_{\eps}}(t))  \, dt \right|^2 \right] \\[0.25cm]
    \leq \, &\left(\sup_{\eps>0} \| \varphi_\eps \|_\infty T_\eps\right)^2 \delta + \int_{|y| \leq M} \rho_\eps(y)^2 \, dy \, \E \left[ \left| \int_0^{T_\eps \wedge \tau_n} \varphi_\eps(\xi_\eps^{x_{\eps}}(t))  \, dt \right|^2 \right].
\end{align*}
Let us analyze the second term in detail. Fubini's theorem, equation \eqref{eq:double_integral_claim_1}, and the triangle inequality enable us to estimate
\begin{align*}
    &\int_{|y| \leq M} \rho_\eps(y)^2 \, dy \, \E \left[ \left| \int_0^{T_\eps \wedge \tau_n} \varphi_\eps(\xi_\eps^{x_{\eps}}(t))  \, dt \right|^2 \right] \\[0.25cm]
    =\, &2 \int_{|y| \leq M} \rho_\eps(y)^2 \, dy \, \E \int_0^{T_\eps \wedge \tau_n} \varphi_\eps(\xi_\eps^{x_{\eps}}(t)) \int_{t}^{T_\eps \wedge \tau_n} \, \varphi_\eps(\xi_\eps^{x_{\eps}}(s)) \, ds \, dt \\[0.25cm]
    =\, &2 \, \E \int_{|y| \leq M} \rho_\eps(y)^2 \int_\R \varphi_\eps(z) \rho_\eps(z)^2 \int_0^{T_\eps \wedge \tau_n} \varphi_\eps(\xi_\eps^{x_{\eps}}(t)) \left[ |\xi_\eps^{x_\eps}(T_\eps \wedge \tau_n) - z| - |\xi_\eps^{x_\eps}(t) - z| \right] dt \, dz \, dy \\[0.25cm]
    \leq\, &2 \, \E \int_{|y| \leq M} \rho_\eps(y)^2 \int_\R \varphi_\eps(z) \rho_\eps(z)^2 \int_0^{T_\eps \wedge \tau_n} \varphi_\eps(\xi_\eps^{x_{\eps}}(t)) \left[ |\xi_\eps^{x_\eps}(T_\eps \wedge \tau_n) - y| - |\xi_\eps^{x_\eps}(t) - y| \right] dt \, dz \, dy \\[0.25cm]
    &+ 4 \, \E \int_{|y| \leq M} \rho_\eps(y)^2 \int_\R \varphi_\eps(z) \rho_\eps(z)^2 \int_0^{T_\eps \wedge \tau_n} \varphi_\eps(\xi_\eps^{x_{\eps}}(t)) |y - z| dt \, dz \, dy
\end{align*}
For the first term of the last inequality we have by equation \eqref{eq:double_integral_claim_2}
\begin{align*}
    &\E \int_{|y| \leq M} \rho_\eps(y)^2 \int_\R \varphi_\eps(z) \rho_\eps(z)^2 \int_0^{T_\eps \wedge \tau_n} \varphi_\eps(\xi_\eps^{x_{\eps}}(t)) \left[ |\xi_\eps^{x_\eps}(T_\eps \wedge \tau_n) - y| - |\xi_\eps^{x_\eps}(t) - y| \right] dt \, dz \, dy \\[0.25cm]
    =\, &\int_\R \varphi_\eps(z)\rho_\eps(z)^2 \, dz \, \E \int_0^{T_\eps \wedge \tau_n} \varphi_\eps(\xi_\eps^{x_{\eps}}(t)) \int_t^{T_\eps \wedge \tau_n} \mathds{1}_{[-M, M]}(\xi_\eps^{x_{\eps}}(s)) \, ds \, dt \\[0.25cm]
    \leq\, &\int_\R \varphi_\eps(z)\rho_\eps(z)^2 \, dz \, \E \int_0^{T_\eps \wedge \tau_n} \varphi_\eps(\xi_\eps^{x_{\eps}}(t)) ( T_\eps \wedge \tau_n - t ) \, dt 
\end{align*}
For the second term it holds by the assumptions
\begin{align*}
    &\hspace{-0.15cm}\left| \E \int_{|y| \leq M} \rho_\eps(y)^2 \int_\R \varphi_\eps(z) \rho_\eps(z)^2 \int_0^{T_\eps \wedge \tau_n} \varphi_\eps(\xi_\eps^{x_{\eps}}(t)) |y - z| dt \, dz \, dy \right| \\[0.25cm]
    \leq\, &\E \int_0^{T_\eps \wedge \tau_n} \varphi_\eps(\xi_\eps^{x_{\eps}}(t)) \, dt \int_{|y| \leq M} \rho_\eps(y)^2 \int_{|z| \leq R} \varphi_\eps(z) \rho_\eps(z)^2 |y-z| \, dz \, dy \\[0.25cm]
    \leq\, &\left(\sup_{\eps>0} \| \varphi_\eps \|_\infty D_\rho\right)^2 (M+R) T_\eps.
\end{align*}
Therefore, combining all these estimates yields
\begin{align*}
    \int_\R \rho_\eps(y)^2 \, dy \, \E \left[ \left| \int_0^{T_\eps \wedge \tau_n} \varphi_\eps(\xi_\eps^{x_{\eps}}(t))  \, dt \right|^2 \right]
    \leq \, &\left(\sup_{\eps>0} \| \varphi_\eps \|_\infty T_\eps\right)^2 \delta + 4\left(\sup_{\eps>0} \| \varphi_\eps \|_\infty D_\rho\right)^2 (M+R) T_\eps \\[0.25cm]
    &+ 2\int_\R \varphi_\eps(z)\rho_\eps(z)^2 \, dz \, \E \int_0^{T_\eps \wedge \tau_n} \varphi_\eps(\xi_\eps^{x_{\eps}}(t)) ( T_\eps \wedge \tau_n - t ) \, dt.
\end{align*}
Letting finally $n \rightarrow \infty$ and using the dominated convergence theorem gives the estimate
\begin{align*}
    \int_\R \rho_\eps(y)^2 \, dy \, \E \left[ \left| \int_0^{T_\eps} \varphi_\eps(\xi_\eps^{x_{\eps}}(t))  \, dt \right|^2 \right]
    \leq \, &\left(\sup_{\eps>0} \| \varphi_\eps \|_\infty T_\eps\right)^2 \delta + 4\left(\sup_{\eps>0} \| \varphi_\eps \|_\infty D_\rho\right)^2 (M+R) T_\eps \\[0.25cm]
    &+ 2 \int_\R \varphi_\eps(z)\rho_\eps(z)^2 \, dz \, \E \int_0^{T_\eps} \varphi_\eps(\xi_\eps^{x_{\eps}}(t)) ( T_\eps - t ) \, dt.
\end{align*}
Hence,
\begin{align*}
    &\E \left[ \left| \frac{1}{T_\eps} \int_0^{T_\eps} \varphi_\eps(\xi_\eps^{x_{\eps}}(t))  \, dt \right|^2 \right] - \left( \frac{1}{C_{\rho_\eps}} \int_\R \varphi_\eps(z)\rho_\eps(z)^2 \, dz \right)^2 \\[0.25cm]
    \leq \, &\frac{1}{C_{\rho_\eps}} \int_\R \varphi_\eps(z)\rho_\eps(z)^2 \, dz \left[ \frac{2}{T_\eps^2} \E \int_0^{T_\eps} \varphi_\eps(\xi_\eps^{x_{\eps}}(t)) ( T_\eps - t ) \, dt - \frac{1}{C_{\rho_\eps}} \int_\R \varphi_\eps(z)\rho_\eps(z)^2 \, dz \right] \\[0.25cm]
    &+ \frac{\delta}{d_{\rho}} \left(\sup_{\eps>0} \| \varphi_\eps \|_\infty \right)^2 + \frac{4}{C_{\rho_\eps} T_\eps} \left(\sup_{\eps>0} \| \varphi_\eps \|_\infty D_\rho\right)^2 (M+R)
\end{align*}
The RHS converges by Lemma \ref{lemma:semigroup_ergodicity} to $\left(\sup_{\eps>0} \| \varphi_\eps \|_\infty \right)^2 \delta/d_\rho$ as $\epstozero$, so that
\begin{equation}
    \limsup_{\epstozero} \left( \E \left[ \left| \frac{1}{T_\eps} \int_0^{T_\eps} \varphi_\eps(\xi_\eps^{x_{\eps}}(t))  \, dt \right|^2 \right] - \left( \frac{1}{C_{\rho_\eps}} \int_\R \varphi_\eps(z)\rho_\eps(z)^2 \, dz \right)^2 \right) \leq  \frac{\delta}{d_\rho} \left(\sup_{\eps>0} \| \varphi_\eps \|_\infty \right)^2.
\end{equation}
Letting $\delta \rightarrow 0$ eventually yields the claim \eqref{eq:mean_ergodic_theorem_limsup_claim}. Thus, the convergence result \eqref{eq:eps_mean_ergodic_theorem_transf} is established for Borel-measurable, nonnegative functions $\varphi_\eps$ with $\sup_{\eps > 0} \| \varphi_\eps \|_{\infty} < \infty$ and $\cup_{\eps > 0} \; \mathrm{supp}(\varphi_\eps) \subset [-R, R]$ for some $R>0$. We extend the convergence result now to measurable, nonnegative functions $\varphi_\eps$ with $\sup_{\eps > 0} \| \varphi_\eps \|_{\infty} < \infty$. For this, let $\delta>0$ be arbitrary and choose again $M>0$ such that for all $\eps>0$
\begin{equation*}
    \int_{|y|>M} \rho_\eps(y)^2 \, dy < \delta.
\end{equation*}
Set $\psi_\eps^{(m)} := \varphi_\eps \mathds{1}_{[-m^*, m^*]}$ with $m^* := m^*(m) := m \vee M$, $m \in \N$. Then we can upper bound
\begin{align*}
    \E \left| \frac{1}{T_\eps} \int_0^{T_\eps} \varphi_\eps(\xi_\eps^{x_{\eps}}(t))  \, dt - \frac{1}{C_{\rho_\eps}} \int_\R \varphi_\eps(y) \rho_\eps(y)^2 \, dy \right|^2
    \leq \, &8 \Bigg[ \E \left| \frac{1}{T_\eps} \int_0^{T_\eps} \varphi_\eps(\xi_\eps^{x_{\eps}}(t))  \, dt - \frac{1}{T_\eps} \int_0^{T_\eps} \psi_\eps^{(m)}(\xi_\eps^{x_{\eps}}(t))  \, dt \right|^2 \\[0.25cm]
    \, &+ \E \left| \frac{1}{T_\eps} \int_0^{T_\eps} \psi_\eps^{(m)}(\xi_\eps^{x_{\eps}}(t))  \, dt - \frac{1}{C_{\rho_\eps}} \int_\R \varphi_\eps(y) \rho_\eps(y)^2 \, dy \right|^2 \\[0.25cm]
    \, &+ \left| \frac{1}{C_{\rho_\eps}} \int_\R \left( \psi_\eps^{(m)}(y) - \varphi_\eps(y) \right) \rho_\eps(y)^2 \, dy \right|^2 \Bigg]
\end{align*}
We estimate the first term as follows
\begin{align*}
    &\E \left| \frac{1}{T_\eps} \int_0^{T_\eps} \varphi_\eps(\xi_\eps^{x_{\eps}}(t))  \, dt - \frac{1}{T_\eps} \int_0^{T_\eps} \psi_\eps^{(m)}(\xi_\eps^{x_{\eps}}(t))  \, dt \right|^2  \\[0.25cm]
    = \, &\E \left| \frac{1}{T_\eps} \int_0^{T_\eps} \varphi_\eps(\xi_\eps^{x_{\eps}}(t)) \mathds{1}_{[-m^*, m^*]^c}(\xi_\eps^{x_{\eps}}(t))  \, dt \right|^2   \\[0.25cm]
    \leq \, & \left(\sup_{\eps>0} \| \varphi_\eps \|_\infty \right)^2 \E \left| 1 - \frac{1}{T_\eps} \int_0^{T_\eps} \mathds{1}_{[-m^*, m^*]}(\xi_\eps^{x_{\eps}}(t)) \, dt \right|^2 \\[0.25cm]
    \leq \, &4\left(\sup_{\eps>0} \| \varphi_\eps \|_\infty \right)^2 \E \left| \frac{1}{T_\eps} \int_0^{T_\eps} \mathds{1}_{[-m^*, m^*]}(\xi_\eps^{x_{\eps}}(t)) \, dt - \frac{1}{C_{\rho_\eps}} \int_{|y| \leq m^*} \rho_\eps(y)^2 \, dy \right|^2 \\[0.25cm]
    &+ 4\left(\sup_{\eps>0} \| \varphi_\eps \|_\infty \right)^2 \left( \frac{1}{C_{\rho_\eps}} \int_{|y| > m^*} \rho_\eps(y)^2 \, dy \right)^2 \\[0.25cm]
    \leq \, &4\left(\sup_{\eps>0} \| \varphi_\eps \|_\infty \right)^2 \left[ \E \left| \frac{1}{T_\eps} \int_0^{T_\eps} \mathds{1}_{[-m^*, m^*]}(\xi_\eps^{x_{\eps}}(t)) \, dt - \frac{1}{C_{\rho_\eps}} \int_{|y| \leq m^*} \rho_\eps(y)^2 \, dy \right|^2 + \left( \frac{\delta}{d_\rho} \right)^2 \right],
\end{align*}
where the last upper bound converges to $4\left(\sup_{\eps>0} \| \varphi_\eps \|_\infty \delta/d_\rho \right)^2$ as $\epstozero$. By our previously established result in this proof, we also have the convergence of the second term for any $m \in \N$ as $\epstozero$. The third term satisfies
\begin{align*}
    \left| \frac{1}{C_{\rho_\eps}} \int_\R \left( \psi_\eps^{(m)}(y) - \varphi_\eps(y) \right) \rho_\eps(y)^2 \, dy \right|^2 
    = \left| \frac{1}{C_{\rho_\eps}} \int_\R \varphi_\eps(y) \mathds{1}_{[-m^*, m^*]^c} \rho_\eps(y)^2 \, dy \right|^2 
    \leq \left( \sup_{\eps>0} \| \varphi_\eps \|_\infty \frac{\delta}{d_\rho} \right)^2.
\end{align*}
In summary, we have
\begin{equation*}
    \limsup_{\epstozero} \E \left| \frac{1}{T_\eps} \int_0^{T_\eps} \varphi_\eps(\xi_\eps^{x_{\eps}}(t))  \, dt - \frac{1}{C_{\rho_\eps}} \int_\R \varphi_\eps(y) \rho_\eps(y)^2 \, dy \right|^2 \leq 40 \left( \sup_{\eps>0} \| \varphi_\eps \|_\infty \frac{\delta}{d_\rho} \right)^2,
\end{equation*}
which implies, by letting $\delta \rightarrow 0$,
\begin{equation*}
     \limsup_{\epstozero} \E \left| \frac{1}{T_\eps} \int_0^{T_\eps} \varphi_\eps(\xi_\eps^{x_{\eps}}(t))  \, dt - \frac{1}{C_{\rho_\eps}} \int_\R \varphi_\eps(y) \rho_\eps(y)^2 \, dy \right|^2 \leq 0.
\end{equation*}
Hence, the convergence result holds for a uniformly bounded sequence $(\varphi_\eps)_{\eps>0}$ of Borel-measurable, nonnegative functions on $\R$, too. The asserted full result is now straightforward and can be accomplished by splitting the test functions into positive and negative part and applying the preceding result to these nonnegative functions.
\end{proof}

\subsection{Convergence in probability result}      \label{subsec:prob_convergence}
This section is concerned with the convergence in probability to zero of the quantity $X_\eps(T_\eps)/\sqrt{T_\eps}$ as $\epstozero$. The final result is contained in Corollary \ref{cor:prob_conv_result} at the end of this section. Most of the crucial proof ideas of this section originate from \cite[Chapter 4]{KMN:2012}, but, as in the preceding section, we must practice caution due to the $\eps$-dependence. 

The first major step towards Corollary \ref{cor:prob_conv_result} is establishing Schauder interior estimates for the parabolic equation
\begin{equation}    \label{eq:parabolic_problem}
    \partial_t u_\eps(t, x) - \mathcal{A}_\eps u_\eps(t, x) = 0, \quad t \in (0, T), \; x \in U,
\end{equation}
where $\mathcal{A}_\eps$ was defined in \eqref{eq:generators}, $T>0$, and $U$ is an open, bounded interval in $\R$. Let us introduce additional terminology that is commonly used in the analysis of parabolic partial differential equations. For $\alpha \in (0,1)$ we say that a function $f \colon (0, T) \times U \to \R$ is locally Hölder-continuous on $(0, T) \times U$ with exponent $\alpha \in (0,1)$ if for any closed subset $D \subset (0, T) \times U$ there exists a constant $H>0$ such that for all $(t, x), (r, y) \in D$ it holds
\begin{equation}
    |f(t,x) -f(r,y)| \leq H \left[|x-y|^2 + |t-r|\right]^{\alpha/2}.
\end{equation}
Furthermore, we define the function $d \colon [0, \infty) \times \overline{U} \to \R$ by
\begin{equation}    \label{eq:distance_function}
    d(t, x) := \text{dist}(x, \partial U) \wedge \sqrt{t}, \quad t \in[0, \infty), \; x \in \overline{U}.
\end{equation}
At last, for $m \in \N_0$ and $\alpha \in (0,1)$ we define the norms
\begin{align}
    \begin{aligned}
    |d^{\, m}, f|_\alpha := &\sup_{(t,x) \in (0, T) \times U} | d(t,x)^m f(t, x)| \\[0.25cm]
    &+ \sup_{\substack{(t, x), (r, y) \in (0, T) \times U \\ (t, x) \neq (r, y)}} | d(t,x) \wedge d(r, y)|^{m+\alpha} \frac{|f(t,x) -f(r,y)|}{\left[|x-y|^2 + |t-r|\right]^{\alpha/2}}.
    \end{aligned}
\end{align}
Observe that whenever the norm $|d^{\, m}, f|_\alpha$ is finite for some $m \in \N_0$ and $\alpha \in (0,1)$, then $f$ is locally Hölder-continuous on $(0, T) \times U$ with exponent $\alpha$.

The following result, whose proof is standard in the literature for partial differential equations (PDEs), see \cite[Chapter 4]{F:1964}, requires, although very technical, merely an explicit tracking of the constant appearing in the upper bound with respect to $\eps$. The acquired estimates on the constant are far from being optimal in any sense, but given our assumptions, this will serve our purposes, as we will soon see.
\begin{Proposition}     \label{prop:schauder_estimates}
Let $\sigma$ be bounded away from zero on $U$ and let the coefficients $\sigma, b_\eps$ of the operator $\mathcal{A}_\eps$ satisfy for an exponent $\alpha \in (0,1)$
\begin{equation}    \label{eq:coefficient_hölder_bounds}
    |d^{\, 0}, \sigma |_\alpha \leq K_\sigma, \quad |d^{\, 1}, b_\eps |_\alpha \leq K_b \, \eps^{-2},
\end{equation}
with some constants $K_\sigma, K_b > 0$. Assume that $u_\eps \in C^{1, 2}((0, T) \times U)$ is a solution of \eqref{eq:parabolic_problem} such that $\sup_{(0, T) \times U} |u_\eps| < \infty$ and $u_\eps$, $\partial_x u_\eps$, $\partial_x^2 u_\eps$, $\partial_t u_\eps$ are locally Hölder-continuous on $(0, T) \times U$ with exponent $\alpha \in (0,1)$. Then there exists a constant $K>0$, only depending on $K_\sigma, K_b$ and $\alpha$, such that for all $\eta > 0$
\begin{equation}    \label{eq:schauder_estimates}
    |d^{\, 0}, u_\eps|_\alpha + |d^{\, 1}, \partial_x u_\eps|_\alpha + |d^{\, 2}, \partial^2_x u_\eps|_\alpha + |d^{\, 2}, \partial_t u_\eps|_\alpha \leq K \eps^{-(10+\eta)} \sup_{(0, T) \times U} |u_\eps|.
\end{equation}
\end{Proposition}
\begin{proof}
    Carefully tracking the constants in \cite[Chapter 4, Section 4]{F:1964} with respect to $\eps$ yields the claim.
\end{proof}

Our aim is to apply Proposition \ref{prop:schauder_estimates} to the Feller transition probability function $P_\eps$ of $X_\eps$, which is given by the relation
\begin{equation}    \label{eq:feller_tpf}
    P_\eps(x, t, A) := \Pb(X_\eps^x(t) \in A), \quad x \in \R, \; t > 0, \; A \in \mathscr{B}(\R). 
\end{equation}
Evidently, it is stochastically continuous, which is to say, by definition, that for all $x \in \R$ and $\delta > 0$
\begin{equation*}
    P_\eps(x, t, B_\delta(x)) = \Pb(|X_\eps^x(t) - x| > \delta) \rightarrow 1, \quad t \rightarrow 0^+.
\end{equation*}
Hence, \cite[Lemma 3.1]{KMN:2012} implies that $(t, x) \mapsto P_\eps(x, t, A)$ is a $\mathscr{B}([0, \infty) \times \R)$-measurable function for any $A \in \mathscr{B}(\R)$, which makes it possible to write down certain integrals. 

Before we proceed, we introduce the following parabolic PDE with Dirichlet initial and boundary conditions and coefficients that do not depend on time
\begin{align}   \label{eq:parabolic_problem_dirichlet}
    \begin{cases}
        \begin{alignedat}{3}
            &\partial_t u_\eps(t, x) - \mathcal{A}_\eps u_\eps(t, x) = 0,     &\quad& t \in (0, T], \; x \in U,          \\[0.25cm]
	    &u_\eps(0, x) = g(x)                                        ,     &\quad& x \in \overline{U},                \\[0.25cm]
            &u_\eps(t, x) = h(t, x)                                     ,     &\quad& t \in (0, T], \; x \in \partial U, \\[0.25cm]
        \end{alignedat}
    \end{cases}
\end{align}
with given functions $g \in C(\overline{U})$ and $h \in C((0, T] \times \partial U)$. If a solution to this initial-boundary value problem exists, then it satisfies the Feynman-Kac formula, cf. \cite[Lemma 3.3, Remark 3.13]{KMN:2012}, that is
\begin{equation}    \label{eq:feynman_kac}
    u_\eps(t, x) = \E \left[ g(X_\eps^x(t)) \mathds{1}_{\{ t \leq \tau_U \}} + h(t-\tau_U, X_\eps^x(\tau_U)) \mathds{1}_{\{ t > \tau_U \}} \right], \quad t \in [0, T], \; x \in \overline{U},
\end{equation}
where $\tau_U := \inf\{t \geq 0 \, | \, X_\eps^x(t) \notin \overline{U} \}$ is the first exit time of the process $X_\eps^{x}$ from $\overline{U}$. This representation can be directly proved with Itô's formula; the needed techniques to carry out this proof can be found in \cite[Chapter 4]{KS:1996}. With these facts at our disposal, we may now state and prove the following lemma.
\begin{Lemma}   \label{lemma:P_eps_PDE_solution}
Assume \ref{ass:assumptions_C} ii). Then, for any fixed $A \in \mathscr{B}(\R)$, the function $(t, x) \mapsto P_\eps(x, t, A)$ is a solution of
\begin{equation}    \label{eq:P_eps_parabolic_problem}
    \partial_t u_\eps(t, x) - \mathcal{A}_\eps u_\eps(t, x) = 0, \quad (t,x) \in D,
\end{equation}
where $D \subset (0, T) \times U$ is an arbitrary subdomain with closure in $(0, T] \times \overline{U}$.
Moreover, for $(t_0, T) \times C \subset (0, T) \times U$ with $t_0 \in (0, T)$ and closed $C$ we have with arbitrary $\eta>0$
\begin{equation}    \label{eq:P_eps_derivative_estimate}
   \sup_{(t, x) \in (t_0, T) \times C} | \partial_x P_\eps(x, t, A) | \leq \frac{K \eps^{-(10+\eta)}}{\mathrm{dist}(\partial C, \partial U) \wedge \sqrt{t_0}}
\end{equation}
where $K>0$ is the same constant as in \eqref{eq:coefficient_hölder_bounds}.
\end{Lemma}
\begin{proof}
    First note that \ref{ass:assumptions_C} ii) implies that the coefficients $\sigma, b_\eps$ of the operator $\mathcal{A}_\eps$ satisfy \eqref{eq:coefficient_hölder_bounds} for any exponent $\alpha \in (0,1)$. Indeed, since the coefficients do not depend on time, it is easy to prove
    \begin{align*}
        |d^{\, 0}, \sigma |_\alpha \leq &\sup_{x \in \overline{U}} | \sigma(x)| + \mathrm{diam}(U)^{\alpha} \sup_{\substack{(x, y) \in \times \overline{U} \\ x \neq y}} \frac{|\sigma(x) - \sigma(y)|}{|x-y|^\alpha}, \\[0.25cm]
        |d^{\, 1}, b_\eps |_\alpha \leq & \mathrm{diam}(U) \sup_{x \in \overline{U}} | b_\eps(x)| + \mathrm{diam}(U)^{1+\alpha} \sup_{\substack{(x, y) \in \times \overline{U} \\ x \neq y}} \frac{|b_\eps(x) - b_\eps(y)|}{|x-y|^\alpha},
    \end{align*}
    which implies the inequalities in \eqref{eq:coefficient_hölder_bounds} after applying \ref{ass:assumptions_C} ii). Observe also that these constants do not depend on $T$.
    With these estimates and the strict nondegeneracy of $\sigma$, the existence of a solution $u_\eps$ to \eqref{eq:parabolic_problem_dirichlet} with locally Hölder-continuous derivatives of all relevant orders for any exponent $\alpha \in (0,1)$ is secured, see \cite[Chapter 3]{F:1964}, whenever the initial and boundary data are continuous.
    We can write for $t \in [0, T]$ and $ x \in \overline{U}$
    \begin{align}  \label{eq:P_eps_decomposition}
        P_\eps(x, t, A) 
        &= \E \left[ \mathds{1}_{A}(X_\eps^{x}(t)) \mathds{1}_{\{ t \leq \tau_U \}} \right] + \Pb\left( X_\eps^{x}(t) \in A, t > \tau_U \right). 
    \end{align}
    The second term can be expressed using the strong Markov property as follows
    \begin{align}
        \begin{aligned}
            \Pb\left( X_\eps^{x}(t) \in A, t > \tau_U \right)
            &= \int_{\left\{ \tau_U < t \right\}} \Pb\left( X_\eps^{x}(t) \in A \, \big| \, \mathcal{F}_{\tau_U} \right)(\omega) \, d\Pb(\omega) \\[0.25cm]
            &= \int_{\Omega} \mathds{1}_{\left\{ \tau_U(\omega) < t \right\}} P_\eps \left( X_\eps^{x}(\tau_U), t - \tau_U, A \right)(\omega) \, d\Pb(\omega) \\[0.25cm]
            &= \int_{\partial U} \int_0^{t} P_\eps\left( y, t - s, A \right) \, d\Pb^{\left(\tau_U, X_\eps^{x}(\tau_U)\right)}(s, y).
        \end{aligned}
    \end{align}
    In the last step we used the $\mathscr{B}([0, \infty) \times \R)$-measurability of $(t, x) \mapsto P_\eps(x, t, A)$. Hence, 
    \begin{equation}    \label{eq:P_eps_second_summand}
        \Pb\left( X_\eps^{x}(t) \in A, t > \tau_U \right) = \E \left[ P_\eps(X_\eps^x(\tau_U), t-\tau_U, A) \mathds{1}_{\{ t > \tau_U \}} \right].
    \end{equation}
    Inserting \eqref{eq:P_eps_second_summand} into \eqref{eq:P_eps_decomposition} gives for $t \in [0, T]$ and $ x \in \overline{U}$
    \begin{equation}    \label{eq:P_eps_feynman_kac}
        P_\eps(x, t, A) = \E \left[ \mathds{1}_{A}(X_\eps^{x}(t)) \mathds{1}_{\{ t \leq \tau_U \}} + P_\eps(X_\eps^x(\tau_U), t-\tau_U, A) \mathds{1}_{\{ t > \tau_U \}} \right],
    \end{equation}
    which is \eqref{eq:feynman_kac}, but with initial and boundary conditions which are only Borel-measurable. We can take care of this with an approximation argument. For this purpose, let $B \in \mathscr{B}([0, T] \times \overline{U})$ be arbitrary. We take two sequences of bounded, piecewise linear functions $k_n^-, k_n^+ \in C([0, \infty) \times \R)$, $n \in \N$, such that
    \begin{equation}    \label{eq:indicator_bounds}
        k_n^- \leq \mathds{1}_B \leq k_n^+, \quad \text{on } [0, \infty) \times \R, \; n \in \N,
    \end{equation}
    and, both, $k_n^-$ and $k_n^+$ converge pointwise from below and from above, respectively, to $\mathds{1}_B$ as $n \rightarrow \infty$. By the discussion at the beginning of the proof we have two sequences of solutions of \eqref{eq:parabolic_problem_dirichlet} $u_{n}^{-}$ and $u_{n}^{+}$, corresponding to $k_n^-$ and $k_n^+$, respectively, with representations given by the Feynman-Kac formula, namely
    \begin{equation}    \label{eq:approx_sequence_feynman_kac}
        u_{n}^{\pm}(t, x) := u_{n, \eps}^{\pm}(t, x) = \E \left[ k_n^{\pm}(\tau_U \wedge t, X_\eps^x(\tau_U \wedge t)) \right], \quad t \in [0, T], \; x \in \overline{U}.
    \end{equation}
    Note here that the initial and boundary functions are summarized by a single function for ease of notation. A well-known result from the PDE literature, e.g., \cite[p.80, Theorem 15]{F:1964}, provides us with subsequences $u_{n_j}^{\pm}$ that converge pointwise on $D$ to some functions $u_\eps^{\pm}$ as $j \rightarrow \infty$, where $D \subset (0, T) \times U$ is an arbitrary subdomain with closure in $(0, T] \times \overline{U}$. Furthermore, these limit functions satisfy $\partial_t u_\eps^\pm - \mathcal{A}_\eps u_\eps^\pm= 0$ on $D$ and have locally Hölder-continuous derivatives of all relevant orders for any exponent $\alpha \in (0,1)$. We will now show that $u_\eps^- = u_\eps^+$ on $D$. By \eqref{eq:indicator_bounds} and \eqref{eq:approx_sequence_feynman_kac} it holds
    \begin{align*}
        \E \left[ k_{n_j}^{-}(\tau_U \wedge t, X_\eps^x(\tau_U \wedge t)) \right] &= u_{n_j}^-(t, x) \\[0.25cm]
        &\leq \E \left[ \mathds{1}_B(\tau_U \wedge t, X_\eps^x(\tau_U \wedge t)) \right] \\[0.25cm]
        &\leq u_{n_j}^+(t, x) = \E \left[ k_{n_j}^{+}(\tau_U \wedge t, X_\eps^x(\tau_U \wedge t)) \right]
    \end{align*}
    for $(t, x) \in [0, T] \times \overline{U}$. Now, on the one hand, we have the pointwise convergence
    \begin{equation*}
        \lim_{j \rightarrow \infty} u_{n_j}^\pm(t, x) = u_\eps^\pm(t, x), \quad (t, x) \in D, 
    \end{equation*}
    and, on the other hand, monotone convergence gives
    \begin{equation*}
        \lim_{j \rightarrow \infty} \E \left[ k_{n_j}^{\pm}(\tau_U \wedge t, X_\eps^x(\tau_U \wedge t)) \right] = \E \left[ \mathds{1}_B(\tau_U \wedge t, X_\eps^x(\tau_U \wedge t)) \right], \quad (t, x) \in [0, T] \times \overline{U}.
    \end{equation*}
    Eventually, this implies
    \begin{equation}
        u_\eps^-(t, x) = u_\eps^+(t, x) = \E \left[ \mathds{1}_B(\tau_U \wedge t, X_\eps^x(\tau_U \wedge t)) \right], \quad (t, x) \in D.
    \end{equation}
    We can extend the preceding result to simple functions and then to $\mathscr{B}([0, T] \times \overline{U})$-measurable functions. Hence, the function $(t, x) \mapsto P_\eps(x, t, A)$ fulfills \eqref{eq:P_eps_parabolic_problem} and, in particular, has locally Hölder-continuous derivatives of all relevant orders for any exponent $\alpha \in (0,1)$. At last, we can derive \eqref{eq:P_eps_derivative_estimate} which is a consequence of \eqref{eq:schauder_estimates}, since
    \begin{equation*}
        \sup_{(t, x) \in (t_0, T) \times C} | \partial_x P_\eps(x, t, A) | \leq \frac{|d^{\, 1}, \partial_x P_\eps( \cdot, \cdot, A) |_\alpha}{\text{dist}(\partial C, \partial U) \wedge \sqrt{t_0}} \leq \frac{K \eps^{-(10+\eta)}}{\text{dist}(\partial C, \partial U) \wedge \sqrt{t_0}}.   \tag*{\qedhere}
    \end{equation*}
\end{proof}
For the remainder of this subsection and the next subsection, we will need the following additional assumption.
\begin{assumptioncustom}{(CLT)}   \label{ass:assumptions_CLT}
There exist numbers $S, \gamma > 0$ such that for all $|y|>S$ it holds
\begin{equation}
    \sign(y) \frac{b(y)}{\overline{\sigma}(y)^2} \leq -\gamma.
\end{equation}
\end{assumptioncustom}
The inequality indicates that the drift coefficient of the limit SDE \eqref{eq:original_sde}, when sufficiently far away from the origin, stays below the threshold $-\gamma$. Note that this assumption already implies Assumptions \ref{ass:assumptions_MET} i) and ii), so that, by Remark \ref{rem:existence_invariant_density}, the invariant densities $\mu$ and $\mu_\eps$ of $X$ and $X_\eps$, respectively, exist. Furthermore, if \ref{ass:assumptions_MET} iii) holds, then it is easy to prove that for all $\eps > 0$
\begin{equation}    \label{eq:bounds_invariant_densities}
    c_\mu \mu \leq \mu_\eps \leq C_\mu \mu \quad \text{on } \R,
\end{equation}
with some constants $c_\mu, C_\mu > 0$ that are independent of $\eps$.

In the proof of the next Proposition \ref{prop:stoch_convergence_khasminskij} we will make use of recurrence properties of the process $X_\eps^{x_0}$. First note that, under Assumptions \ref{ass:assumptions_C} and \ref{ass:assumptions_MET}, the process $X_\eps$ is recurrent for each $\eps> 0$, i.e., for any $y \in \R$
\begin{equation}
    \Pb(\tau^{y}\left(X_\eps^{x_0}\right) < \infty) = 1,
\end{equation}
where we recall $\tau^{y}\left(X_\eps^{x_0}\right) = \inf \{t>0 \, | \, X_\eps^{x_0}(t) = y\}$. This is a consequence of \cite[Lemma 3.9, Remark 15]{KMN:2012}, cf. with \cite[Example 3.10]{KMN:2012}, as well. We even have positive recurrence of $X_\eps$ uniformly in $\eps$. To be precise, for any $y \in \R$ it holds
\begin{equation}    \label{eq:positive_recurrence}
    \sup_{\eps > 0} \E \tau^{y}\left(X_\eps^{x_0}\right) < \infty.
\end{equation}
Indeed, by Lemma \ref{lemma:random_time_lemma}, we have for all $\eps > 0$ and $y \in \R$
\begin{align*}
    \E \tau^{y}\left(X_\eps^{x_0}\right) = \E \tau^{f_\eps(y)}(\xi_\eps^{f_\eps(x_0)}) \leq 2 D_\rho |f_\eps(x_0)-f_\eps(y)| \leq 2 D_\rho ( |f(x_0)|+|f(y)| ).
\end{align*}
\begin{Proposition} \label{prop:stoch_convergence_khasminskij}
    If Assumptions \ref{ass:assumptions_C}, \ref{ass:assumptions_MET} and \ref{ass:assumptions_CLT} hold, then for any initial condition $x_0 \in \R$
    \begin{equation}
        \Pb\left( \left| X_\eps^{x_0}(T_\eps) \right| > R(\eps) \right) \rightarrow 0, \quad \epstozero,
    \end{equation}
    where $R(\eps) > 0$ is a sequence with $R(\eps) \rightarrow \infty$ as $\epstozero$ which will be determined in the proof.
\end{Proposition}
\begin{proof}
    Fix $I := (-y, y)$ for some $y > 1$. Consider the initial datum $x_0 \in I$ first, and let $T > 0$ be arbitrary. Due to Lemma \ref{lemma:P_eps_PDE_solution} it holds for all $x \in B_{r(\eps)}(x_0)$, where $r(\eps) \downarrow 0$ as $\epstozero$ will be determined soon, and $t_0 \in (0, T)$
    \begin{equation}
        \left| P_\eps(x, T, B_{R(\eps)}(0)^c) - P_\eps(x_0, T, B_{R(\eps)}(0)^c) \right| \leq \frac{K \eps^{-(10+\eta)}}{\mathrm{dist}(\partial B_{r(\eps)}(x_0), \partial I) \wedge \sqrt{t_0}} |x - x_0|.
    \end{equation}
    Evidently, as $\epstozero$ it holds $\text{dist}(\partial B_{r(\eps)}(x_0), \partial I) \rightarrow \mathrm{dist}(x_0, \partial I) = |x_0 -y| > 0$, so that choosing $r(\eps) := \eps^{11+\eta}$ gives
    \begin{equation}    \label{eq:P_eps_mean_value}
        \left| P_\eps(x, T, B_{R(\eps)}(0)^c) - P_\eps(x_0, T, B_{R(\eps)}(0)^c) \right| \leq K \eps,
    \end{equation}
    with a new constant $K>0$ independent of $\eps$. As the next step towards the proof we may estimate, using \eqref{eq:bounds_invariant_densities} and stationarity as in \eqref{eq:stationarity},
    \begin{equation}
        \int_\R P_\eps(x, T, B_{R(\eps)}(0)^c) \mu_\eps(x) \, dx = \int_{B_{R(\eps)}(0)^c} \mu_\eps(x) \, dx \leq C_\mu \int_{B_{R(\eps)}(0)^c} \mu(x) \, dx,
    \end{equation}
    and using \eqref{eq:P_eps_mean_value} it follows
    \begin{align*}
        \int_\R P_\eps(x, T, B_{R(\eps)}(0)^c) \mu_\eps(x) \, dx &\geq c_\mu \int_{B_{r(\eps)}(x_0)} P_\eps(x, T, B_{R(\eps)}(0)^c) \mu(x) \, dx \\[0.25cm]
        &\geq c_\mu \int_{B_{r(\eps)}(x_0)} \mu(x) \, dx \left[ P_\eps(x_0, T, B_{R(\eps)}(0)^c) - K \eps \right].
    \end{align*}
    Rearranging terms yields
    \begin{equation}    \label{eq:prior_P_eps_estimate}
        \Pb\left( \left| X_\eps^{x_0}(T) \right| > R(\eps) \right) = P_\eps(x_0, T, B_{R(\eps)}(0)^c) \leq K \eps + \frac{C_\mu \int_{B_{R(\eps)}(0)^c} \mu(x) \, dx}{c_\mu \int_{B_{r(\eps)}(x_0)} \mu(x) \, dx},
    \end{equation}
    for all $T > 0$ and $x_0 \in I$. Observe that, by continuity, we have on the one hand 
    \begin{equation}    \label{eq:inner_ball_estimate}
        \lim_{\epstozero} \frac{1}{2 r(\eps)} \int_{B_{r(\eps)}(x_0)} \mu(x) \, dx = \mu(x_0) > 0, \quad \epstozero.
    \end{equation}
    On the other hand, by Assumption \ref{ass:assumptions_CLT}, we have for $x > S$
    \begin{align*}
        \mu(x) = \frac{1}{Z} \exp \left( 2\int_0^x \frac{b(y)}{\overline \sigma(y)^2} \, dy \right) \leq C(b, \overline \sigma, S) \exp \left( 2\int_S^x \frac{b(y)}{\overline \sigma(y)^2} \, dy \right) \leq C(b, \overline \sigma, S, \gamma) \exp(-2 \gamma x),
    \end{align*}
    and, similarly, for $x < -S$ but with different sign. For sufficiently small $\eps > 0$ so that $R(\eps) \geq S$, the last estimate implies
    \begin{equation}    \label{eq:outer_ball_estimate}
        \int_{B_{R(\eps)}(0)^c} \mu(x) \, dx \leq C(b, \overline \sigma, S, \gamma) \exp(-2 \gamma R(\eps)).
    \end{equation}
    Combining \eqref{eq:prior_P_eps_estimate}, \eqref{eq:inner_ball_estimate}, and \eqref{eq:outer_ball_estimate} eventually yields
    \begin{equation}    \label{eq:prob_convergence_estimate}
        \Pb\left( \left| X_\eps^{x_0}(T) \right| > R(\eps) \right) \leq K \eps + C(b, \overline \sigma, S, \gamma, \mu) \frac{\exp(-2\gamma R(\eps))}{r(\eps) \mu(x_0)} \rightarrow 0, \quad \epstozero,
    \end{equation}
    when choosing, say, $R(\eps) = \mathcal{O}(\eps^{-\eta})$ with $\eta > 0$. This preceding result also holds when $T = T_\eps$, since the appearing estimates did not depend on $T$. \\[0.25cm]
    Now consider $x_0 \in I^c$, assume w.l.o.g $x_0 > y$, and let $\delta > 0$. By \eqref{eq:positive_recurrence} and Markov's inequality we can choose $T_0 := T_0(\delta, x_0, y) > 0$, independent of $\eps$, such that for all $\eps > 0$
    \begin{equation}    \label{rev:markov_estimate}
         \Pb(\tau^{y-1}\left(X_\eps^{x_0}\right) > T_0) \leq \frac{\sup_{\eps > 0} \E \tau^{y-1}\left(X_\eps^{x_0}\right)}{T_0}< \delta.
    \end{equation}
    Notice that we use $y-1 \in I$ here, because we soon want to utilize \eqref{eq:prob_convergence_estimate} which was only proved for interior points $x_0 \in I$. Next, we have for sufficiently small $\eps > 0$ such that $T_\eps \geq T_0$
    \begin{align*}
        \Pb\left( \left| X_\eps^{x_0}(T_\eps) \right| > R(\eps) \right) 
        &\leq \delta + \Pb\left( \left| X_\eps^{x_0}(T_\eps) \right| > R(\eps), \tau^{y-1}\left(X_\eps^{x_0}\right) \leq T_\eps \right). 
    \end{align*}
    The last term can be estimated using the strong Markov property again and \eqref{eq:prob_convergence_estimate} as follows
    \begin{align*}
        \Pb\left( \left| X_\eps^{x_0}(T_\eps) \right| > R(\eps), \tau^{y-1}\left(X_\eps^{x_0}\right) \leq T_\eps \right)
        &= \int_0^{T_\eps} P_\eps\left( y-1, T_\eps - s, B_{R(\eps)}(0)^c \right) \, d\Pb^{\tau^{y-1}\left(X_\eps^{x_0}\right)}(s) \\[0.25cm]
        &\leq K \eps + C(b, \overline \sigma, S, \gamma, \mu) \frac{\exp(-2 \gamma R(\eps))}{r(\eps) \mu(y-1)}.
    \end{align*}
    Note in the step before that $\Pb(\tau^{y-1}\left(X_\eps^{x_0}\right) = 0) = 0$ since $x_0 > y$. Hence, as $\epstozero$
    \begin{equation*}
        \lim_{\epstozero} \Pb\left( \left| X_\eps^{x_0}(T_\eps) \right| > R(\eps) \right) \leq \delta,
    \end{equation*}
    which yields the claim for the case $x_0 \in I^c$, as well.
\end{proof}
\begin{Remark}
    In the proof of Proposition \ref{prop:stoch_convergence_khasminskij}, instead of $R(\eps) = \mathcal{O}(\eps^{-\eta})$ one can also choose $R(\eps) = -(2\gamma)^{-1} \mathcal{O}(\log(r(\eps) \eps^\eta))$ for $\epstozero$, $\eta > 0$, and get the same convergence result.
\end{Remark}
\begin{Corollary}   \label{cor:prob_conv_result}
    Let Assumptions \ref{ass:assumptions_C}, \ref{ass:assumptions_MET} and \ref{ass:assumptions_CLT} hold and let $T_\eps = \mathcal{O}(\eps^{-\eta})$ as $\epstozero$ with $\eta > 0$. Then for all $\delta > 0$ and any $x_0 \in \R$
    \begin{equation}
        \lim_{\epstozero} \Pb\left( \frac{\left| X_\eps^{x_0}(T_\eps) \right|}{\sqrt{T_\eps}} > \delta \right) = 0.
    \end{equation}
    In other words, $ X_\eps^{x_0}(T_\eps)/\sqrt{T_\eps}$ converges in probability to zero as $\epstozero$.
\end{Corollary}

\subsection{\texorpdfstring{An $\eps$}{Lg}-dependent central limit theorem}     \label{subsec:central_limit_theorem}
In this section, we want to extend the results of the preceding sections by outlining and proving an "$\eps$-dependent" CLT under certain additional assumptions, i.e., for a given function $h_\eps \colon \R \to \R$ we want to prove 
\begin{equation}    \label{eq:eps_clt_intro}
    \frac{1}{\sqrt{T_\eps}} \int_0^{T_\eps} h_\eps(X_\eps^{x_0}(t)) \, dt \weakconv \mathcal{N}(0, \tau^2), \quad \epstozero,
\end{equation}
where $\weakconv$ denotes weak convergence of measures in $\R$ and $\tau^2>0$ is the asymptotic variance. The idea is based on the technique of the Poisson equation, which is classical in the literature, cf.\ \cite{KLO:2012} or \cite{MST:2010}, and goes as follows.

Fix a continuous function $h_\eps \in L^1(\mu_\eps)$ satisfying $\int_\R h_\eps(x) \mu_\eps(x) \, dx = 0$ with $\mu_\eps$ as in Remark \ref{rem:existence_invariant_density}. If there exists a solution $\Phi_\eps \in \mathcal{D}(\mathcal{A}_\eps) \subset C^2(\R)$ to the Poisson equation
\begin{equation}    \label{eq:eps_poisson_equation}
    - \mathcal{A}_\eps \Phi_\eps = h_\eps,
\end{equation}
where $\mathcal{A}_\eps$ is the differential operator introduced in \eqref{eq:generators} with domain $\mathcal D(\mathcal A_\eps)$, then applying the Itô formula to $\Phi_\eps$ yields
\begin{equation}    \label{eq:clt_poisson_method}
    \frac{1}{\sqrt{T_\eps}} \int_0^{T_\eps} h_\eps(X_\eps^{x_0}(t)) \, dt = \frac{\Phi_\eps(x_0) - \Phi_\eps(X_\eps^{x_0}(T_\eps))}{\sqrt{T_\eps}} + \frac{1}{\sqrt{T_\eps}} \int_0^{T_\eps} \sigma(X_\eps^{x_0}(t)) \Phi'_\eps(X_\eps^{x_0}(t)) \, dW(t).
\end{equation}
Inspecting equation \eqref{eq:clt_poisson_method} we notice that, in order to obtain the CLT \eqref{eq:eps_clt_intro}, we need to control the first term in such a way that it vanishes in probability as $\epstozero$ and use a CLT for the remaining stochastic integral. The latter claim is contained in \cite[Exercise (IV.3.33)]{RY:2013} so that we merely state the result in the way that we require.
\begin{Proposition} \label{prop:eps_clt_revuz_yor}
For every $\eps > 0$ let $\phi_\eps \colon \R \to \R$ be a measurable function such that $\int_0^{T_\eps} \phi_\eps(X_\eps^{x_0}(t))^2 \, dt < \infty$, $\Pb$-a.s. If there exists $\tau > 0$ such that
\begin{equation}
    \frac{1}{T_\eps} \int_0^{T_\eps} \phi_\eps(X_\eps^{x_0}(t))^2 \, dt \probconv \tau^2, \quad \epstozero,
\end{equation}
then
\begin{equation}
    \frac{1}{\sqrt{T_\eps}} \int_0^{T_\eps} \phi_\eps(X_\eps^{x_0}(t)) \, dW(t) \weakconv \mathcal{N}(0, \tau^2), \quad \epstozero.
\end{equation}
\end{Proposition}

To apply this proposition to \eqref{eq:clt_poisson_method} we first provide a lemma concerning the properties of $\Phi_\eps$ and its derivative. The main ideas for this lemma are borrowed from the proof of Lemma 1.17 in \cite{K:2004} and are modified to our setting. 
\begin{Lemma}   \label{lemma:properties_poisson_solution}
Let Assumptions \ref{ass:assumptions_C}, \ref{ass:assumptions_MET}, and \ref{ass:assumptions_CLT} hold. Assume that the functions $h_\eps \in L^1(\mu_\eps)$ are continuous and satisfy 
\begin{align}   \label{eq:h_eps_assumption}
    \int_\R h_\eps(x) \mu_\eps(x) \, dx = 0, \quad  \sup_{\eps > 0} \|h_\eps\|_\infty < \infty. 
\end{align}
Then the sequence of functions $\Phi'_\eps$, where $\Phi_\eps$ solves \eqref{eq:eps_poisson_equation}, is uniformly bounded.
\end{Lemma}
\begin{proof}
In this particular one-dimensional setting the function $\Phi_\eps$ is given by
\begin{equation} \label{eq:poisson_eps}
    \Phi_\eps(x) = - \int_0^x \frac{2}{\sigma(y)^2 \mu_\eps(y)} \int_{-\infty}^y h_\eps(z) \mu_\eps(z) \, dz \, dy, \quad x \in \R,
\end{equation}
Its first derivative is
\begin{equation} \label{eq:derivative_poisson}
    \Phi'_\eps(x) = - \frac{2}{\sigma(x)^2 \mu_\eps(x)} \int_{-\infty}^x h_\eps(z) \mu_\eps(z) \, dz, \quad x \in \R.
\end{equation}
Recall the inequality \eqref{eq:bounds_invariant_densities}
\begin{equation*}
    c_\mu \mu \leq \mu_\eps \leq C_\mu \mu \quad \text{on } \R,
\end{equation*}
and observe that by the integral condition in \eqref{eq:h_eps_assumption}
\begin{align*}
    \Phi'_\eps(x) = -\frac{2}{\sigma(x)^2 \mu_\eps(x)} \int_{-\infty}^x h_\eps(z) \mu_\eps(z) \, dz 
    = \frac{2}{\sigma(x)^2 \mu_\eps(x)} \int_{x}^\infty h_\eps(z) \mu_\eps(z) dz, \quad x \in \R.
\end{align*}
Now, let $S, \gamma$ be as in \ref{ass:assumptions_CLT}. By using \ref{ass:assumptions_C} ii), the inequality in \ref{ass:assumptions_CLT}, and \eqref{eq:h_eps_assumption}, we obtain for $x > S$
\begin{align*}
    |\Phi'_\eps(x)| &\leq \frac{2 C_\mu}{c_\mu} \int_{x}^\infty \frac{|h_\eps(z)|}{\sigma(z)^2} \exp \left( \int_x^z \frac{2 b(y)}{\overline{\sigma}(y)^2} \, dy \right) dz \\[0.25cm]
    &\leq \frac{2 C_\mu \sup_{\eps > 0} \|h_\eps\|_\infty}{c_\mu a^2} \int_{x}^\infty \, \exp \left(- 2\gamma (z-x) \right) dz
    = \frac{C_\mu \sup_{\eps > 0} \|h_\eps\|_\infty}{\gamma c_\mu a^2}.
\end{align*}
A similar estimate is also true for $x < -S$. For $x \in [-S, S]$ it obviously holds
\begin{equation*}
    |\Phi'_\eps(x)| \leq \frac{\sup_{\eps > 0} \|h_\eps\|_\infty}{c_\mu a^2 \mu(x)},
\end{equation*}
which implies the claim.
\end{proof}

For the next result, we also fix a continuous function $h \in L^1(\mu)$ with $\int_\R h(x) \mu(x) \, dx = 0$, and consider
\begin{equation} \label{eq:poisson}
    \Phi(x) = - \int_0^x \frac{2}{\overline{\sigma}(y)^2 \mu(y)} \int_{-\infty}^y h(z) \mu(z) \, dz \, dy, \quad x \in \R,
\end{equation}
which solves the Poisson equation $- \mathcal{A} \Phi = h$ with the differential operator $\mathcal{A}$ given in \eqref{eq:generators}.
\begin{Theorem} \label{theo:eps_clt}
Let Assumptions \ref{ass:assumptions_C}, \ref{ass:assumptions_MET}, and \ref{ass:assumptions_CLT}, and $T_\eps = \mathcal{O}(\eps^{-\eta})$ as $\epstozero$ with $\eta > 0$ hold. Assume that the functions $h_\eps \in L^1(\mu_\eps)$ are continuous and satisfy 
\begin{align}
    \int_\R h_\eps(x) \mu_\eps(x) \, dx = 0, \quad  \sup_{\eps > 0} \|h_\eps\|_\infty < \infty. 
\end{align} 
Additionally, assume that
\begin{equation}    \label{eq:asymptotic_variance_convergence}
    \tau_\eps^2 := \int_\R \sigma(x)^2 \Phi'_\eps(x)^2 \mu_\eps(x) \, dx \rightarrow \int_\R \overline{\sigma}(x)^2 \Phi'(x)^2 \mu(x) \, dx =: \tau^2, \quad \epstozero.
\end{equation}
Then
\begin{equation}
    \frac{1}{\sqrt{T_\eps}} \int_0^{T_\eps} h_\eps(X_\eps^{x_0}(t)) \, dt \weakconv \mathcal{N}(0, \tau^2), \quad \epstozero\,.
\end{equation}
\end{Theorem}
\begin{proof}
Recall the equation \eqref{eq:clt_poisson_method}
\begin{equation*}
    \frac{1}{\sqrt{T_\eps}} \int_0^{T_\eps} h_\eps(X_\eps^{x_0}(t)) \, dt = \frac{\Phi_\eps(x_0) - \Phi_\eps(X_\eps^{x_0}(T_\eps))}{\sqrt{T_\eps}} + \frac{1}{\sqrt{T_\eps}} \int_0^{T_\eps} \sigma(X_\eps^{x_0}(t)) \Phi'_\eps(X_\eps^{x_0}(t)) dW(t).
\end{equation*}
By virtue of Corollary \ref{cor:prob_conv_result}, the first term on the right-hand side vanishes in probability as $\epstozero$ because $\Phi'_\eps$ is uniformly bounded by Lemma \ref{lemma:properties_poisson_solution}. For the second term, we first bound
\begin{equation*}
    \E \left| \frac{1}{T_\eps} \int_0^{T_\eps} \sigma(X_\eps^{x_0}(t))^2 \Phi'_\eps(X_\eps^{x_0}(t))^2 dt - \tau^2 \right| \leq \E \left| \frac{1}{T_\eps} \int_0^{T_\eps} \sigma(X_\eps^{x_0}(t))^2 \Phi'_\eps(X_\eps^{x_0}(t))^2 dt - \tau_\eps^2 \right| + | \tau_\eps^2 - \tau^2 |,
\end{equation*}
and then apply Theorem \ref{theo:eps_mean_ergodic_theorem} to get
\begin{equation*}
    \E \left| \frac{1}{T_\eps} \int_0^{T_\eps} \sigma(X_\eps^{x_0}(t))^2 \Phi'_\eps(X_\eps^{x_0}(t))^2 dt - \tau_\eps^2 \right| \rightarrow 0, \quad \epstozero.
\end{equation*}
By condition \eqref{eq:asymptotic_variance_convergence} it therefore follows
\begin{equation*}
    \E \left| \frac{1}{T_\eps} \int_0^{T_\eps} \sigma(X_\eps^{x_0}(t))^2 \Phi'_\eps(X_\eps^{x_0}(t))^2 dt - \tau^2 \right| \rightarrow 0,  \quad \epstozero.
\end{equation*}
Hence, the claim is established in view of Proposition \ref{prop:eps_clt_revuz_yor} and Slutzky's Lemma.
\end{proof}

\section{Application to parameter estimation}   \label{sec:parameter_estimation}
This section is devoted to a simple application of the previous limit theorems to a statistical parameter estimation problem. Consider the following overdamped Langevin diffusion in one dimension with a linear and oscillatory term in the drift and constants $\alpha, \sigma > 0$
\begin{equation} \label{eq:multiscale_langevin_linear_exa}
    dX_\eps(t) = \left[ -\alpha X_\eps(t) - \frac{1}{\eps} \sin\left( \frac{X_\eps(t)}{\eps} \right) \right] dt  + \sqrt{2 \sigma} dW(t), \quad X_\eps(0) = x_0.
\end{equation}%
We omit the superscript $x_0$ in $X_\eps^{x_0}$ for this section. Using homogenization theory, see \cite[Chapter 3]{BLP:1978} or \cite{PS:2007}, it can be shown that the law $\Pb^{X_\eps}$ converges weakly to $\Pb^{X}$ in $C([0,T];\R^d)$ as $\epstozero$ for fixed $T>0$, where $X$ is the solution of the SDE
\begin{equation} \label{eq:limit_langevin_linear_exa}
    dX(t) = -\alpha K X(t) dt  + \sqrt{2 \sigma K} d\overline{W}(t), \quad X(0) = x_0.
\end{equation}
Here, the constant factor $K>0$ emerges from the cell problem of the homogenization and equals
\begin{equation}
    K := \frac{1}{Z^+ Z^-}, \quad Z^{\pm} := \frac{1}{2\pi} \int_0^{2 \pi} \exp \left( \pm \frac{\cos(y)}{\sigma} \right) \, dy.
\end{equation}
It is fairly easy to see that all the conditions in \ref{ass:assumptions_C}, \ref{ass:assumptions_MET}, and \ref{ass:assumptions_CLT} are satisfied.
In particular, the invariant densities $\mu_\eps$ and $\mu$ of $X_\eps$ and $X$, respectively, exist and are given by
\begin{align*}    \label{eq:eps_invariant_density_langevin}
    \mu_\eps(x) = \frac{1}{Z_\eps} \exp \left( -\frac{\alpha}{2 \sigma} x^2 + \frac{1}{\sigma} \cos\left( \frac{x}{\eps} \right) \right), \quad
    \mu(x) = \frac{1}{Z} \exp \left(  -\frac{\alpha}{2 \sigma} x^2 \right), \quad x \in \R,
\end{align*}
where $Z$ and $Z_\eps$ are normalization constants.

We want to estimate the parameter $\vartheta := \alpha K$ appearing in \eqref{eq:limit_langevin_linear_exa} with observations coming in the form of a single trajectory from \eqref{eq:multiscale_langevin_linear_exa}. For simplicity, we assume that we know the value $\overline{\sigma} := \sigma K$, e.g., from a prior estimation procedure. Such diffusion parameter estimation problems with given multiscale data were analyzed, for example, in \cite{KPK:2013, MP:2018, AGPSZ:2021}. We propose the following simple minimum distance estimator based on the characteristic function of the invariant density $\mu$ at the point $1$
\begin{equation}    \label{eq:mde_exa_def}
    \hvt_{T_\eps}(X_\eps) := \arginf_{\vt \in \Theta} \left| \frac{1}{T_\eps} \int_0^{T_\eps} \exp(i X_\eps(t)) \, dt - \exp\left( -\frac{\overline{\sigma}}{2\vartheta} \right) \right|.
\end{equation}
This estimator is similar in spirit to the newly proposed estimator in \cite{BKP:2025} where the authors give a more comprehensive analysis and review of the parameter estimation problem at hand. Using differentiation, it is straightforward to obtain an explicit expression for this estimator, namely
\begin{equation}    \label{eq:mde_exa_closed_form}
    \hvt_{T_\eps}(X_\eps) = - \frac{\overline{\sigma}}{2 \log\left( \frac{1}{T_\eps} \int_0^{T_\eps} \exp(i X_\eps(t)) \, dt \right)}.
\end{equation}
Consider the Poisson equations
\begin{equation}    \label{eq:poisson_exa}
    -\mathcal{A}_\eps \Phi_\eps = h_\eps, \quad -\mathcal{A} \Phi = h, \quad \text{on } \R,
\end{equation}
with the differential operators $\mathcal{A}_\eps$ and $\mathcal{A}$ corresponding to \eqref{eq:multiscale_langevin_linear_exa} and \eqref{eq:limit_langevin_linear_exa}, respectively, and with the functions
\begin{align*}
    h_\eps(x) := \exp(ix) - \int_\R \exp(iy) \mu_\eps(y) \, dy,  \quad h(x) := \exp(ix) - \int_\R \exp(iy) \mu(y) \, dy, \quad x \in \R.
\end{align*}
We need the following lemma.
\begin{Lemma}   \label{lemma:char_function_convergence}
Let $\alpha, \sigma, R > 0$. Then
\begin{equation}    \label{eq:char_functions_estimate}
    \| h_\eps - h \|_\infty \leq C(\alpha, \sigma) \exp\left( -\frac{2\sigma}{\alpha \eps^2} \right), 
\end{equation}
and
\begin{equation}    \label{eq:Phi_estimate}
    \sup_{x \in [-R, R]} |\Phi_\eps(x) - \Phi(x)| \rightarrow 0, \quad \text{as } \epstozero.
\end{equation}
\end{Lemma}
\begin{proof}
In this proof, the appearing constants $C(\alpha, \sigma)$ and $C(\alpha, \sigma, R)$ may change from line to line, but they will, in any case, stay independent of $\eps$. We will first prove the exponential convergence in \eqref{eq:char_functions_estimate}. First, notice that we can write for any $z \in \R$
\begin{align*}
    h_\eps(z) - h(z) 
    &= \int_{\R} \exp(ix) \left[ \mu_\eps(x) - \mu(x) \right] dx \\[0.25cm]
    &= \int_{\R} \exp(ix) \left[ \mu_\eps(x) - \frac{Z^+}{Z_\eps} \exp\left(-\frac{\alpha}{2\sigma} x^2\right) \right] dx + \int_{\R} \exp(ix) \left[ \frac{Z^+}{Z_\eps} \exp\left(-\frac{\alpha}{2\sigma} x^2\right) - \mu(x) \right] dx \\[0.25cm]
    &= \frac{1}{Z_\eps} \int_{\R} \exp(ix) \exp\left(-\frac{\alpha}{2\sigma} x^2\right) \left[ \exp\left( \frac{\cos(x/\eps)}{\sigma} \right) - Z^+ \right] dx \\[0.25cm]
    &\hspace{0.5cm} + \frac{Z Z^+ - Z_\eps}{Z Z_\eps} \int_{\R^d} \exp(ix) \exp\left(-\frac{\alpha}{2\sigma} x^2\right) \, dx.
\end{align*}
Observe that
\begin{equation*}
    Z_\eps - Z Z^+ = \int_{\R^d} \exp\left(-\frac{\alpha}{2\sigma} x^2\right) \left[ \exp\left( \frac{\cos(x/\eps)}{\sigma} \right) - \frac{1}{2\pi} \int_0^{2\pi} \exp\left(\frac{\cos(y)}{\sigma}\right) \, dy \right] dx.
\end{equation*}
We will prove that
\begin{equation}    \label{eq:constants_conv}
    |Z_\eps - Z Z^+| \leq C(\alpha, \sigma) \exp\left( -\frac{2\sigma}{\alpha \eps^2} \right).
\end{equation}
For this, we can use a series expansion and the binomial theorem to calculate the integral
\begin{align*}
    \frac{1}{2 \pi} \int_0^{2 \pi} \exp \left( \frac{\cos(y)}{\sigma} \right) dy
    &= \frac{1}{2\pi} \sum_{n=0}^\infty \frac{1}{\sigma^n n!} \int_0^{2 \pi} \cos(y)^n dy \\[0.25cm]
    &= \frac{1}{2\pi} \sum_{n=0}^\infty \frac{1}{(2 \sigma)^n n!} \int_0^{2 \pi} \left( \exp \left( iy \right) + \exp \left( -iy \right) \right)^n dy \\[0.25cm]
    &= \frac{1}{2\pi} \sum_{n=0}^\infty \sum_{k=0}^n \frac{1}{(2 \sigma)^n n!} \binom{n}{k} \int_0^{2 \pi} \exp \left( i(n-2k)y \right) dy \\[0.25cm]
    &= \frac{1}{2\pi} \sum_{m=0}^\infty \sum_{k=0}^{2m} \frac{1}{(2 \sigma)^{2m}} \frac{1}{k! (2m -k)!} \int_0^{2 \pi} \exp \left( 2i(m-k)y \right) dy = \sum_{m=0}^\infty \frac{1}{(2 \sigma)^{2m} (m!)^2}.
\end{align*}
Similarly for $\eps > 0$ and $x \in \R$
\begin{align*}
    \exp \left( \frac{\cos(x/\eps)}{\sigma} \right) 
    &= \sum_{n=0}^\infty \sum_{k=0}^n \frac{1}{(2 \sigma)^{n} n!} \binom{n}{k} \exp \left( i(n-2k)x/\eps \right) \\[0.25cm]
    &= \sum_{m=0}^\infty \frac{1}{(2 \sigma)^{2m} (m!)^2} + \sum_{m=0}^\infty \sum_{\substack{ k=0 \\ k \neq m}}^{2m} \frac{\exp \left( 2i(m-k)x/\eps  \right)}{(2 \sigma)^{2m} k! (2m -k)!} \\[0.25cm]
    &\hspace{0.35cm} + \sum_{m=0}^\infty \sum_{k=0}^{2m+1} \frac{\exp \left( i(2(m-k)+1)x/\eps \right)}{(2 \sigma)^{2m+1} k! (2m-k+1)!}.
\end{align*}
We will estimate the last two terms. For the middle term in the last appearing equation we have with Fubini's theorem, the fact that $\E \exp(itZ) = \exp(-t^2/2)$ when $Z \sim \mathcal{N}(0,1)$, $|m-k| \geq 1$. and the binomial theorem again
\begin{align}
\begin{aligned} \label{eq:middle_term_estimate}
    &\left| \int_\R \exp\left(-\frac{\alpha}{2\sigma} x^2\right) \sum_{m=0}^\infty \sum_{\substack{ k=0 \\ k \neq m}}^{2m} \frac{\exp \left( 2i(m-k)x/\eps \right)}{(2 \sigma)^{2m} k! (2m -k)!} \, dx \right| \\[0.25cm]
    &\leq \sum_{m=0}^\infty \sum_{\substack{ k=0 \\ k \neq m}}^{2m} \frac{1}{(2 \sigma)^{2m} k! (2m -k)!} \left| \int_\R \exp\left(-\frac{\alpha}{2\sigma} x^2\right) \exp \left( 2i(m-k)x/\eps \right) \, dx \right| \\[0.25cm]
    &\leq C(\alpha, \sigma) \sum_{m=0}^\infty \sum_{\substack{ k=0 \\ k \neq m}}^{2m} \frac{1}{(2 \sigma)^{2m} k! (2m -k)!} \exp\left( -\frac{2(m-k)^2 \sigma}{\alpha \eps^2} \right) \\[0.25cm]
    &\leq C(\alpha, \sigma) \exp\left( -\frac{2\sigma}{\alpha \eps^2} \right) \sum_{m=0}^\infty \frac{1}{\sigma^{2m} (2m)!}.
\end{aligned}
\end{align}
An analog estimate holds for the third term, namely
\begin{align*}
    \left| \int_\R \exp\left(-\frac{\alpha}{2\sigma} x^2\right) \sum_{m=0}^\infty \sum_{k=0}^{2m+1} \frac{\exp \left( i(2(m-k)+1)x/\eps \right)}{(2 \sigma)^{2m+1} k! (2m-k+1)!} \, dx \right| \leq C(\alpha, \sigma) \exp\left( -\frac{2\sigma}{\alpha \eps^2} \right) \sum_{m=0}^\infty \frac{1}{\sigma^{2m+1} (2m+1)!}.
\end{align*}
With these equations and estimates we therefore arrive at \eqref{eq:constants_conv}. Note also that we can get a similar estimate for the term
\begin{align*}
    &\int_{\R} \exp(ix) \exp\left(-\frac{\alpha}{2\sigma} x^2\right) \left[ \exp\left( \frac{\cos(x/\eps)}{\sigma} \right) - Z^+ \right] dx \\[0.25cm]
    &= \int_{\R} \exp(ix) \exp\left(-\frac{\alpha}{2\sigma} x^2\right) \left[ \exp\left( \frac{\cos(x/\eps)}{\sigma} \right) - \frac{1}{2\pi} \int_0^{2\pi} \exp\left(\frac{\cos(y)}{\sigma}\right) \, dy \right] dx
\end{align*}
This is because in \eqref{eq:middle_term_estimate} we only have to replace $\exp(-\alpha x^2/2\sigma)$ with $\exp(ix) \exp(-\alpha x^2/2\sigma)$ so that we end up with the characteristic function of an unnormalized Gaussian again. With these arguments, we can therefore conclude the proof of \eqref{eq:char_functions_estimate}. As a byproduct, we have just proved that $\mu_\eps$ converges weakly to $\mu$ with convergence rates given in terms of their characteristic functions, that is
\begin{equation}
    \left| \int_{\R} \exp(ix) \mu(x) \, dx - \int_{\R} \exp(ix)\mu(x) \, dx \right| \leq C(\alpha, \sigma) \exp\left( -\frac{2\sigma}{\alpha \eps^2} \right).
\end{equation}
Moving on to claim \eqref{eq:Phi_estimate}, we define
\begin{equation}
    H_\eps(y) :=  \int_{-\infty}^y h_\eps(z) \mu_\eps(z) \, dz, \quad H(y) :=  \int_{-\infty}^y h(z) \mu(z) \, dz, \quad y \in \R.
\end{equation}
Then for $y \in \R$
\begin{align*}
    |H_\eps(y) - H(y)| 
    &= \left| \int_{-\infty}^y (h_\eps(z) - h(z)) \mu_\eps(z) \, dz + \int_{-\infty}^y h(z) (\mu_\eps(z) - \mu(z)) \, dz \right| \\[0.25cm]
    &\leq \norm{h_\eps - h}_{\infty} + \left| \int_{-\infty}^y h(z) (\mu_\eps(z) - \mu(z)) \, dz \right|.
\end{align*}
The first term vanishes exponentially fast as we have seen above. The second term vanishes, by weak convergence of $\mu_\eps \lambda^1$ because the function $h \mathds{1}_{(-\infty, y]}$ is $\lambda^1$-a.e. continuous and bounded. Thus, $H_\eps$ converges pointwise to $H$ as $\epstozero$. Now fix $x \in [-R, R]$ and consider $x > 0$. The case $x < 0$ works the same. Recall the solution formulas for $\Phi_\eps$ and $\Phi$ in \eqref{eq:poisson_eps} and \eqref{eq:poisson}. Then it holds
\allowdisplaybreaks
\begingroup
\begin{align*}
    |\Phi_\eps(x) - \Phi(x)| 
    &= \left| \int_0^x \frac{H_\eps(y)}{\sigma \mu_\eps(y)} - \frac{H(y)}{\overline{\sigma} \mu(y)} \, dy \right| \\[0.25cm]
    &= \left| \int_0^x \frac{H_\eps(y) - H(y)}{\sigma \mu_\eps(y)} \, dy + \frac{1}{\sigma} \int_0^x H(y) \left( \frac{1}{\mu_\eps(y)} - \frac{1}{K \mu(y)} \right) \, dy \right| \\[0.25cm]
    &\leq \frac{1}{\sigma c_\mu} \int_0^R \frac{|H_\eps(y) - H(y)|}{\mu(y)} \, dy + \left| \frac{1}{\sigma} \int_0^x H(y) \left( \frac{1}{\mu_\eps(y)} - \frac{1}{K \mu(y)} \right) \, dy \right|.
\end{align*}
\endgroup
Note that $|H_\eps - H| \leq 4$ for all $\eps > 0$, so that the first term converges by dominated convergence. For the second term, we observe that
\begin{align*}
    \int_0^x H(y) \left( \frac{1}{\mu_\eps(y)} - \frac{1}{K \mu(y)} \right) \, dy
    &= Z_\eps \int_0^x H(y) \exp\left(\frac{\alpha}{2\sigma} y^2\right) \left( \exp(-\cos(y/\eps)/\sigma) - Z^- \right) \, dy \\[0.25cm]
    &\hspace{0.35cm} + \int_0^x (Z_\eps Z^- - ZZ^+ Z^-) \exp\left(\frac{\alpha}{2\sigma} y^2\right) \, dy
\end{align*}
Very similar arguments that previously achieved \eqref{eq:constants_conv} can be used to obtain
\begin{equation}    \label{eq:Phi_oscillatory_estimate}
    \int_0^x H(y) \exp\left(\frac{\alpha}{2\sigma} y^2\right) \left( \exp(-\cos(y/\eps)/\sigma) - Z^- \right) \, dy \leq C(\alpha, \sigma, R) \eps.
\end{equation}
Indeed, all we have to do is upper bound the following term in \eqref{eq:middle_term_estimate}
\begin{equation*}
    \left| \int_0^x H(y) \exp\left(\frac{\alpha}{2\sigma} y^2\right) \exp \left( 2i(m-k)y/\eps \right) \, dx \right| \leq C(\alpha, \sigma, R) \eps,
\end{equation*}
which can be accomplished with integration by parts. Hence, by \eqref{eq:constants_conv} and \eqref{eq:Phi_oscillatory_estimate}, we have for $\epstozero$
\begin{equation}
    \sup_{x \in [-R, R]} \int_0^x H(y) \left( \frac{1}{\mu_\eps(y)} - \frac{1}{K \mu(y)} \right) \, dy \rightarrow 0.
\end{equation}
The claim \eqref{eq:Phi_estimate} is finally proved.
\end{proof}
The preceding lemma allows us now to prove asymptotic properties of the minimum distance estimator for $\epstozero$.
\begin{Proposition}     \label{prop:mde_exa_result}
Let $T_\eps = \mathcal{O}(\eps^{-\eta})$ as $\epstozero$ for some $\eta > 0$. Then, under the true parameter $\vt_0 = \alpha_0 K$ with $\alpha_0 > 0$, as $\epstozero$
\begin{align}
    \label{eq:mde_consistency}
    \hvt_{T_\eps}(X_\eps) &\probconv \vt_0, \\[0.25cm]
    \label{eq:mde_asy_normality}
    \sqrt{T_\eps} ( \hvt_{T_\eps}(X_\eps) - \vt_0 ) &\weakconv \mathcal{N}\left(0, \left( \frac{\overline \sigma \tau}{2 \vt_0 \log(\vt_0)^2} \right)^2 \right),
\end{align}
where
\begin{equation}    \label{eq:mde_asymptotic_variance}
    \tau^2 = 2\overline{\sigma} \int_\R \Phi'(x)^2 \mu(x) \, dx.
\end{equation}
\end{Proposition}
\begin{proof}
Fix $\vt_0 = \alpha_0 K$ with $\alpha_0 > 0$ as the true parameter. First observe that
\begin{align*}
    \frac{1}{T_\eps} \int_0^{T_\eps} \exp(iX_\eps(t)) \, dt - \exp\left( -\frac{\overline{\sigma}}{2\vt_0} \right) = \frac{1}{T_\eps} \int_0^{T_\eps} h_\eps(X_\eps(t)) \, dt + h_\eps(x_0) - h(x_0),
\end{align*}
which converges in $L^2\left( \Omega, \mathcal{F}, \Pb \right)$ to zero as $\epstozero$ by Theorem \ref{theo:eps_mean_ergodic_theorem}, Remark \ref{rem:mean_ergodic_theorem_original_quantities}, and Lemma \ref{lemma:char_function_convergence}. Hence,
\begin{equation*}
    \hvt_{T_\eps}(X_\eps) = - \frac{\overline{\sigma}}{2 \log\left( \frac{1}{T_\eps} \int_0^{T_\eps} \exp(i X_\eps(t)) \, dt \right)} \probconv \vt_0.
\end{equation*}
Furthermore, we can write
\begin{equation}    \label{rev:clt_decomposition_example}
    \frac{1}{\sqrt{T_\eps}} \int_0^{T_\eps} \exp(iX_\eps(t)) - \exp\left( -\frac{\overline{\sigma}}{2\vt_0} \right) \, dt = \frac{1}{\sqrt{T_\eps}} \int_0^{T_\eps} h_\eps(X_\eps(t)) \, dt + \sqrt{T_\eps} (h_\eps(x_0) - h(x_0))
\end{equation}
By the exponential convergence of $h_\eps$ to $h$ and the assumption that $T_\eps = \mathcal{O}(\eps^{-\eta})$ as $\epstozero$, we recognize that the second term goes to zero as $\epstozero$. The first term will converge weakly to $\mathcal{N}(0, \tau^2)$ by Theorem \ref{theo:eps_clt} as soon as we establish \eqref{eq:asymptotic_variance_convergence}, that is
\begin{align}   \label{eq:asymptotic_variance_convergence_exa}
    2 \sigma \int_\R \Phi'_\eps(x)^2 \mu_\eps(x) \, dx \rightarrow 2 \overline{\sigma} \int_\R \Phi'(x)^2 \mu(x) \, dx, \quad \epstozero.
\end{align}
Integration by parts and using the Poisson equations in \eqref{eq:poisson_exa} shows that
\begin{align}   \label{eq:dirichlet_form}
    \sigma \int_\R \Phi'_\eps(x)^2 \mu_\eps(x) \, dx = \int_\R \Phi_\eps(x) h_\eps(x) \mu_\eps(x) \, dx, \quad
    \overline{\sigma} \int_\R \Phi'_\eps(x)^2 \mu_\eps(x) \, dx = \int_\R \Phi(x) h(x) \mu(x) \, dx.
\end{align}
We choose a sufficiently large $R>0$ such that the integral
\begin{equation}
    \int_{[-R, R]^c} |\Phi_\eps(x) - \Phi(x)| \,  \mu(x) \, dx 
\end{equation}
gets arbitrarily small uniformly in $\eps$, which is, by Lemma \ref{lemma:properties_poisson_solution}, possible due to the uniform integrability of the integrand. We can perform the following splitting
\begin{align*}
    \int_\R \Phi_\eps(x) h_\eps(x) \mu_\eps(x) \, dx - \int_\R \Phi(x) h(x) \mu(x) \, dx 
    &=  \int_\R \Phi_\eps(x) (h_\eps(x) - h(x)) \mu_\eps(x) \, dx \\[0.25cm]
    &\hspace{0.35cm} + \int_\R (\Phi_\eps(x) - \Phi(x)) h(x) \mu_\eps(x) \, dx \\[0.25cm]
    &\hspace{0.35cm} + \int_\R \Phi(x) h(x) (\mu_\eps(x) - \mu(x))\, dx.
\end{align*}
The first term vanishes by virtue of \eqref{eq:char_functions_estimate} and Lemma \ref{lemma:properties_poisson_solution}. The third term vanishes due to weak convergence of $\mu_\eps \lambda^1$ and uniform integrability of $\Phi h$ with respect to $\mu_\eps \lambda^1$. The second term can be estimated as follows
\begin{align*}
    \left| \int_\R (\Phi_\eps(x) - \Phi(x)) h(x) \mu_\eps(x) \, dx \right| = 2 \left[ \sup_{x \in [-R, R]} |\Phi_\eps(x) - \Phi(x)| + \int_{[-R,R]^c} |\Phi_\eps(x) - \Phi(x)| \, \mu_\eps(x) \, dx \right].
\end{align*}
These terms go to zero as pointed out before and by \eqref{eq:Phi_estimate}. This proves \eqref{eq:asymptotic_variance_convergence_exa}, which, in turn, implies
\begin{equation*}
    \frac{1}{\sqrt{T_\eps}} \int_0^{T_\eps} \exp(iX_\eps(t)) - \exp\left( -\frac{\overline{\sigma}}{2\vt_0} \right) \, dt \weakconv \mathcal{N}(0, \tau^2).
\end{equation*}
The final claim \eqref{eq:mde_asy_normality} is now a simple application of the delta method to the function $g(\vt) := -\overline{\sigma}/(2 \log(\vt))$ and yields
\begin{align*}
    \sqrt{T_\eps} ( \hvt_{T_\eps}(X_\eps) - \vt_0 ) 
    &= \sqrt{T_\eps} \left[ g\left( \frac{1}{T_\eps} \int_0^{T_\eps} \exp(iX_\eps(t)) \, dt \right) - g\left(\exp\left( -\frac{\overline{\sigma}}{2\vt_0} \right) \right) \right] \\[0.25cm]
    &\weakconv \mathcal{N}\left(0, \left( \frac{\overline \sigma \tau}{2 \vt_0 \log(\vt_0)^2} \right)^2 \right), \quad \epstozero. \tag*{\qedhere}
\end{align*}
\end{proof}

\section{Conclusion}\label{sec:conclusion}
In this work, we presented new limit theorems for one-dimensional diffusion processes $X_\eps$ that depend on a small scale parameter $\eps$, for which the processes $X_\eps$ converge weakly to a limit diffusion process $X$. Under sufficient local regularity conditions on the coefficients of the respective SDEs, recurrence assumptions, and sufficient explosion rates of the time horizon $T_\eps \rightarrow \infty$ as $\epstozero$, we were able to prove a mean ergodic theorem and a central limit theorem for time integrals involving bounded test functions. We then applied these new theoretical results to a simple statistical parameter estimation problem under model misspecification in a continuous-time setting. The proposed minimum distance estimator, based on the characteristic function of the limit invariant density, turned out to be consistent and asymptotically normal under observations of the process $X_\eps$ in the coupled long-time and small-scale limiting regime, that is, $T_\eps \rightarrow \infty$ as $\epstozero$.

We see two immediate chances for future research. The first one is the extension of the limit theorems to unbounded test functions that are integrable with respect to the limit invariant measure. This is certainly more useful for statistical estimation purposes where moments of processes are an ubiquitous occurrence. We suspect that, in this case, one must shift away the attention from the mean ergodic theorem and pursue an almost sure ergodic theorem. After all, if the ultimate aim is the consistency or asymptotic normality of estimators, then it is irrelevant if one has the mean ergodic or the almost sure ergodic theorem -- even a stochastic ergodic theorem in the sense of convergence in probability would suffice. 

The second direction for future research is the generalization of the limit theorems to multidimensional processes, e.g., multiscale fast-slow diffusion systems, cf. \cite{BLP:1978, PV:2001, RX:2021}. This poses considerable challenges from several points of view. For example, the techniques for the proof of the mean ergodic theorem are already not extendable to the multidimensional case due to the local time process, whose effective use is limited to the one-dimensional setting. Another difficulty is the question about the convergence of the asymptotic variances in the central limit theorem, which is often proved using reversibility, a rather weak assumption in dimension one, cf. equation \eqref{eq:dirichlet_form}, but fairly restrictive in higher dimensions. When dropping the reversibility assumption, then condition \eqref{eq:asymptotic_variance_convergence} remains rather unclear in more general settings. Thus far we were only able to prove it case by case. Despite these preliminary problems, the convergence in probability results in subsection \ref{subsec:prob_convergence} would be easier to establish since they heavily rely on Schauder interior estimates, which can be proved analogously in several dimensions. Furthermore, we believe that, by combining certain ideas from this material and the newly established estimates in \cite{RX:2021}, it may be possible to obtain a mean ergodic theorem for bounded test functions.

\section*{Acknowledgements}
The authors J.I.B.\ and S.K.\ acknowledge funding by the Deutsche Forschungsgemeinschaft (DFG, German Research Foundation) - Project number 442047500 through the Collaborative Research Center "Sparsity and Singular Structures" (SFB 1481).

\bibliography{v2/ref_v2.bib}

\begin{thebibliography}{10}

\bibitem{AGPSZ:2021}
Assyr Abdulle, Giacomo Garegnani, Grigorios~A. Pavliotis, Andrew~M. Stuart, and Andrea Zanoni.
\newblock Drift {estimation} of {multiscale} {diffusions} {based} on {filtered} {data}.
\newblock {\em Found. Comput. Math.}, 23(1):33--84, February 2021.
\newblock \href {https://doi.org/10.1007/s10208-021-09541-9} {\path{doi:10.1007/s10208-021-09541-9}}.

\bibitem{BLP:1978}
Alain Bensoussan, J.-L. Lions, and George Papanicolaou.
\newblock {\em Asymptotic analysis for periodic structures}, volume~5 of {\em Studies in mathematics and its applications}.
\newblock North-Holland Pub. Co. ; sole distributors for the U.S.A. and Canada, Elsevier North-Holland, Amsterdam, 1st edition, 1978.
\newblock \href {https://doi.org/10.1090/chel/374} {\path{doi:10.1090/chel/374}}.

\bibitem{BKP:2025}
Jaroslav~I. Borodavka, Sebastian Krumscheid, and Grigorios~A. Pavliotis.
\newblock A minimum distance estimator approach for misspecified ergodic processes.
\newblock Submitted to SIAM Multiscale Model. Simul., preprint arxiv:2506.12432, 2025.

\bibitem{DZ:1996}
G.~Da~Prato and J.~Zabczyk.
\newblock {\em Ergodicity for Infinite Dimensional Systems}.
\newblock London Mathematical Society Lecture Note Series. Cambridge University Press, Cambridge, 1996.
\newblock \href {https://doi.org/10.1017/CBO9780511662829} {\path{doi:10.1017/CBO9780511662829}}.

\bibitem{F:1964}
A.~Friedman.
\newblock {\em Partial Differential Equations of Parabolic Type}.
\newblock Prentice-Hall, New Jersey, 1964.
\newblock \href {https://doi.org/10.2307/3613664} {\path{doi:10.2307/3613664}}.

\bibitem{GS:1968}
I.~I. Gikhman and A.V. Skorokhod.
\newblock {\em Stochastic Differential Equations}.
\newblock Springer, New York, 1968.
\newblock Russian version.

\bibitem{KS:1996}
Ioannis Karatzas and Steven~E. Shreve.
\newblock {\em Brownian motion and stochastic calculus}.
\newblock Springer, New York, 1996.
\newblock \href {https://doi.org/10.1007/978-1-4612-0949-2} {\path{doi:10.1007/978-1-4612-0949-2}}.

\bibitem{K:1966b}
R.~Z. Khas'minskii.
\newblock A limit theorem for the solutions of differential equations with random right-hand sides.
\newblock {\em Theory of Probability \& Its Applications}, 11(3):390--406, 1966.
\newblock \href {https://doi.org/10.1137/1111038} {\path{doi:10.1137/1111038}}.

\bibitem{K:1966a}
R.~Z. Khas'minskii.
\newblock On stochastic processes defined by differential equations with a small parameter.
\newblock {\em Theory of Probability \& Its Applications}, 11(2):211--228, 1966.
\newblock \href {https://doi.org/10.1137/1111018} {\path{doi:10.1137/1111018}}.

\bibitem{KMN:2012}
R.~Z. Khas'minskii, G.~N. Milshtein, and M.B. Nevelson.
\newblock {\em Stochastic stability of differential equations}.
\newblock Stochastic modelling and applied probability 66. Springer, Heidelberg, 2nd edition, 2012.
\newblock \href {https://doi.org/10.1007/978-3-642-23280-0} {\path{doi:10.1007/978-3-642-23280-0}}.

\bibitem{K:1949}
A.~I. Khinchin.
\newblock {\em Mathematical Foundations of Statistical Mechanics}.
\newblock Dover Publications, New York, 1949.
\newblock Translated from the Russian by G. Gamow.

\bibitem{KLO:2012}
Thomasz Komorowski, Claudio Landim, and Stefano Olla.
\newblock {\em Fluctuations in Markov processes: time symmetry and martingale approximation}.
\newblock Grundlehren der mathematischen Wissenschaften 345. Springer, Berlin, Heidelberg, 2012.
\newblock \href {https://doi.org/10.1007/978-3-642-29880-6} {\path{doi:10.1007/978-3-642-29880-6}}.

\bibitem{K:2018}
S.~Krumscheid.
\newblock Perturbation-based inference for diffusion processes: {Obtaining} effective models from multiscale data.
\newblock {\em Math. Models Methods Appl. Sci.}, 28(8):1565--1597, 2018.
\newblock \href {https://doi.org/10.1142/S0218202518500434} {\path{doi:10.1142/S0218202518500434}}.

\bibitem{KPK:2013}
S.~Krumscheid, G.~A. Pavliotis, and S.~Kalliadasis.
\newblock Semiparametric {drift} and {diffusion} {estimation} for {multiscale} {diffusions}.
\newblock {\em Multiscale Model. Simul.}, 11(2):442--473, January 2013.
\newblock \href {https://doi.org/10.1137/110854485} {\path{doi:10.1137/110854485}}.

\bibitem{K:2004}
Yury~A. Kutoyants.
\newblock {\em Statistical inference for ergodic diffusion processes}.
\newblock Springer series in statistics. Springer, London, 2004.
\newblock \href {https://doi.org/10.1007/978-1-4471-3866-2} {\path{doi:10.1007/978-1-4471-3866-2}}.

\bibitem{LS:2016}
Tony Lelièvre and Gabriel Stoltz.
\newblock Partial differential equations and stochastic methods in molecular dynamics.
\newblock {\em Acta Numerica}, 25:681–880, 2016.
\newblock \href {https://doi.org/10.1017/S0962492916000039} {\path{doi:10.1017/S0962492916000039}}.

\bibitem{MP:2018}
Theodoros Manikas and Anastasia Papavasiliou.
\newblock Diffusion parameter estimation for the homogenized equation.
\newblock {\em Multiscale Model. Simul.}, 17:675--695, 2018.
\newblock \href {https://doi.org/10.1137/18M120138X} {\path{doi:10.1137/18M120138X}}.

\bibitem{MST:2010}
Jonathan~C. Mattingly, Andrew~M. Stuart, and M.~V. Tretyakov.
\newblock Convergence of numerical time-averaging and stationary measures via {P}oisson equations.
\newblock {\em SIAM J. Numer. Anal.}, 48(2):552--577, 2010.
\newblock \href {https://doi.org/10.1137/090770527} {\path{doi:10.1137/090770527}}.

\bibitem{PSV:1976}
G.~C. Papanicolaou, D.~Stroock, and S.~R.~S. Varadhan.
\newblock Martingale approach to some limit theorems.
\newblock {\em Papers from the Duke Turbulence Conference (Duke Univ., Durham, N.C., 1976)}, Duke Univ. Math. Ser. III:Paper No. 6, ii+120 pp, 1976.
\newblock MR0461684.
\newblock URL: \url{https://mathscinet.ams.org/mathscinet/relay-station?mr=0461684}.

\bibitem{PPS:2009}
A.~Papavasiliou, G.~A. Pavliotis, and A.~M. Stuart.
\newblock Maximum likelihood drift estimation for multiscale diffusions.
\newblock {\em Stoch. Process. Their Appl.}, 119(10):3173--3210, 2009.
\newblock \href {https://doi.org/10.1016/j.spa.2009.05.003} {\path{doi:10.1016/j.spa.2009.05.003}}.

\bibitem{P:1999}
E.~Pardoux.
\newblock Homogenization of linear and semilinear second order parabolic pdes with periodic coefficients: A probabilistic approach.
\newblock {\em J. Funct. Anal.}, 167(2):498--520, 1999.
\newblock \href {https://doi.org/10.1006/jfan.1999.3441} {\path{doi:10.1006/jfan.1999.3441}}.

\bibitem{PV:2001}
E.~Pardoux and A.~Yu. Veretennikov.
\newblock On the {P}oisson equation and diffusion approximation 1.
\newblock {\em Ann. Probab.}, 29(3):1061--–1085, 2001.
\newblock \href {https://doi.org/10.1214/aop/1015345596} {\path{doi:10.1214/aop/1015345596}}.

\bibitem{PS:2007}
G.~A. Pavliotis and A.~M. Stuart.
\newblock Parameter estimation for multiscale diffusions.
\newblock {\em J. Stat. Phys.}, 127:741–--781, 2007.
\newblock \href {https://doi.org/10.1007/s10955-007-9300-6} {\path{doi:10.1007/s10955-007-9300-6}}.

\bibitem{PS:2008}
G.~A. Pavliotis and A.~M. Stuart.
\newblock {\em Multiscale Methods: Averaging and Homogenization}.
\newblock Texts in applied mathematics 53. Springer, New York, 2008.
\newblock \href {https://doi.org/10.1007/978-0-387-73829-1} {\path{doi:10.1007/978-0-387-73829-1}}.

\bibitem{PRZ:2025}
Grigorios~A. Pavliotis, Sebastian Reich, and Andrea Zanoni.
\newblock Filtered data based estimators for stochastic processes driven by colored noise.
\newblock {\em Stoch. Process. Their Appl.}, 181:104558, 2025.
\newblock \href {https://doi.org/10.1016/j.spa.2024.104558} {\path{doi:10.1016/j.spa.2024.104558}}.

\bibitem{RY:2013}
D.~Revuz and M.~Yor.
\newblock {\em Continuous Martingales and Brownian Motion}, volume 293 of {\em Grundlehren der mathematischen Wissenschaften}.
\newblock Springer, Berlin, Heidelberg, 2013.
\newblock \href {https://doi.org/10.1007/978-3-662-06400-9} {\path{doi:10.1007/978-3-662-06400-9}}.

\bibitem{RX:2021}
Michael R{\"o}ckner and Longjie Xie.
\newblock {Diffusion approximation for fully coupled stochastic differential equations}.
\newblock {\em Ann. Probab.}, 49(3):1205 -- 1236, 2021.
\newblock \href {https://doi.org/10.1214/20-AOP1475} {\path{doi:10.1214/20-AOP1475}}.

\bibitem{S:2010}
Tamar Schlick.
\newblock {\em Molecular Modeling and Simulation: An Interdisciplinary Guide}, volume~21 of {\em Interdisciplinary Applied Mathematics}.
\newblock Springer, New York, 2010.
\newblock \href {https://doi.org/10.1007/978-1-4419-6351-2} {\path{doi:10.1007/978-1-4419-6351-2}}.

\bibitem{S:1989}
P.~Sundar.
\newblock Ergodic solutions of stochastic differential equations.
\newblock {\em Stoch. stoch. rep.}, 28(2):65--83, 1988.
\newblock \href {https://doi.org/10.1080/17442508908833584} {\path{doi:10.1080/17442508908833584}}.

\end{thebibliography}

\end{document}